%\magnification=1200
% This file is supposed to imitate the -ras.tex- file which contained
% the macros -insertRaster- and -RasterBox- for including SUN raster files
% in TeX.  But here we assume that the original raster files have been
% converted to PostScript, so we use -psfig- rather than -sunbitmap-.
% There was a bug in the way -ras.tex- scaled these bitmaps, and this file
% attempts to preserve this bug so that the resulting figures are identical
% to the raster figures in the original paper.

\input psfig
%\magnification=1200

\newif\ifboxfigure      % set to true if you want the figure boxed
\boxfigurefalse

\def\BoxIt#1#2{%        % put a box around #1, leaving a gap of #2
	\vbox{\hrule
	\hbox{\vrule\kern#2\vbox{\kern#2#1\kern#2}\kern#2\vrule}
		   \hrule}}

\def\insertRaster #1 pixels #2  by #3 scaled #4 {
%  #1 is the raster file
%  #2 is the width of the picture in pixels
%  #3 is the height of the picture in pixels
%  #4 is scaling 500 means half size 2000 means double size 1000 is actual  
%%%%			\medskip
			 \hbox {%

			 \hss
			 \RasterBox {#1} {#2} {#3} {#4}
			 \hss
			 }%
}

\def\RasterBox #1 #2 #3 #4{

% dimen0 and dimen1 are the dimensions of the hbox which are setup in
% the same way as -ras.tex- unfortunately psfig includes images in a 
% different manner from sunbitmap, so it must be scaled using different
% values, hence dimen2 and dimen3 and MAGIC number 933

\dimen0=#2pt   % #2 is in pixels, so must multiply by 1.1076 when
			  % used
\dimen1=#3pt
\divide\dimen0 by 11076 % really want 1.1076
\multiply\dimen0 by 10 % (because #4 is out of 1000 )
\divide\dimen1 by 11076
\multiply\dimen1 by 10
 \dimen2=#3pt
 \divide\dimen2 by 933
 \dimen3=#2pt
 \divide\dimen3 by 933
\setbox4=\hbox to #4\dimen0{
 \vbox to #4\dimen1{
 \vss
 \psfig{figure=#1,height=#4\dimen2,width=#4\dimen3}
 }
 \hss
 }
 \ifboxfigure\BoxIt{\box4}{0pt}
 \else\box4
 \fi
 }

% \insertRaster towrrab.pic pixels 500 by 500 scaled 400

\parskip=\baselineskip
\baselineskip=\baselineskip 
\hoffset = 0.5 true in
\hsize 5.5  true in
\voffset = 0.5 true in
\vsize 8 true in
\def\endofproof{\hfill $\#$ } 
\def\inchap{\vskip 1.25 true in}
\def\insec{\bigskip\bigskip}
\def\tskip{\vskip -5 true pt}
\vskip -10 true pt
\font \bigfont=cmb10 at 17 true pt
\font \medfont=cmb10 at 15 true pt
\inchap
\centerline{\bigfont On Post Critically Finite Polynomials}
\centerline{\bigfont Part One: Critical Portraits}
\bigskip
\bigskip
\bigskip
\centerline{Alfredo Poirier}
\centerline{Math Department}
\centerline{Suny@StonyBrook}
\centerline{StonyBrook, NY 11794-3651}
\bigskip
\bigskip
\bigskip
{\sl We extend the work of 
{\it Bielefeld, Fisher and Hubbard} 
on Critical Portraits 
(see [BFH] and [F]) to the case of arbitrary postcritically finite polynomials. 
This determines an effective classification of postcritically finite polynomials as dynamical systems.} 
\bigskip
\bigskip
This paper is the first in a series of two 
based on the author's thesis (see [P]) 
which deals with the classification of postcritically finite polynomials. 
In this first part we conclude the study of critical portraits initiated by Fisher (see [F]) 
and continued by Bielefeld, Fisher and Hubbard (see [BHF]). 
As an application of our results, 
we give in the second part of this series 
necessary and sufficient conditions for the realization of Hubbard Trees. 
\vfill\eject
\inchap
\centerline{\bigfont Chapter I}
\centerline{\bigfont Basic Concepts and Main Results}
\inchap 
{In the first three sections of this introductory chapter we define the concept of critically marked polynomials 
and state their main combinatorial properties. 
Our definition extends the concept presented in [F] and [BFH] by including the possibility of periodic critical points. 
This definition differs slightly from that given in the above references in the strictly preperiodic case, 
but our results are the same. 
This small modification will later be useful, 
because some proofs will be simplified. 
\smallskip
Our definition is supported by a number of examples given in Section 4. 
We remark here that the `hierarchic selections' in the construction, are essential only 
to the marking corresponding to Fatou set critical cycles. 
Here they are needed in order to guarantee uniqueness for the polynomial 
with specified critical portrait. 
(Compare Example 4.5, and see the remark following Lemma II.2.4). 
The inclusion of the `hierarchic selection' for 
Julia set critical points was made to uniformize notation and 
is not essential 
(compare [BFH] where all critical points are in the Julia set).}
\insec
\centerline{\medfont 1. Preliminaries.}
\insec 
{\bf 1.1.} 
Let $P$ be a polynomial of degree $d>1$ with $\Omega(P)$ the set of critical points. 
For $M \subset {\bf C}$ denote by 
${\cal O}(M)=\cup^\infty_{n=0} P^{\circ n}(M)$ the {\it orbit of M}. 
If the orbit ${\cal O}(\Omega(P))$ of the critical set is finite, 
we say that $P$ is {\it postcritically finite (PCF)}. 
It follows that every critical point of $P$ is periodic or preperiodic. 
We call the orbit ${\cal O}(\omega)$ of a periodic critical point $\omega$ (if any) a 
{\it critical cycle}. 
In this postcritically finite case a criterion to decide when a preperiodic (or periodic) point is in the Fatou set is as follows.
\smallskip
{\it A preperiodic point is in the Fatou set if and only if it eventually maps to a critical cycle.}
\medskip
If $P$ is postcritically finite, then the Julia set $J(P)$ and the filled in Julia set $K(P)$ of $P$ are connected and locally connected (see [M] Theorem 17.5).
As there are no wandering domains for the Fatou components of this polynomial $P$, 
each bounded Fatou component contains exactly one point $z$ 
(called its {\it center}) which eventually maps to a critical point. 
If we map this component $U(z)$ onto the unit disk by an uniformizing Riemann map 
$\phi$ with $\phi(z)=0$, 
we can talk about {\it internal rays} in $U(z)$ defined as the preimages of radial segments under $\phi$. 
Because we are in the locally connected case those internal rays can be extended up to the boundary.
\medskip
In the case of the basin of attraction of $\infty$, 
if the polynomial is monic and centered, 
the uniformizing Riemann map can be chosen tangent to the identity at $\infty$. 
These rays are called {\it external rays}, 
and satisfy the condition $P(R_\theta)=R_{d\theta}$. 
\smallskip
In general, let $\omega \mapsto P(\omega) \mapsto \dots \mapsto P^{\circ n}(\omega)=\omega$ be a critical cycle. 
Then $P^{\circ n}:U(\omega) \mapsto U(\omega)$ is a degree ${\cal D}>1$ cover of itself 
(${\cal D}$ is the product of the local degree of elements in the orbit ${\cal O}(\omega)$, 
and $U(\omega)$ the Fatou component  with center $\omega$). 
It follows then that the uniformizing Riemann map $\phi_\omega$ can be chosen so that 
$$
\phi_\omega(z)^{\cal D}=\phi_\omega(P^{\circ n}(z)).
$$
In this case the Riemann map is known as the {\it B\"ottcher coordinate} (compare [M] Theorem 6.7). 
This coordinate is uniquely defined up to conjugation with a $({\cal D}-1)^{th}$ root of unity. 
In particular, it is easy to see that there are exactly ${\cal D}-1$ {\it `fixed' internal rays}, 
i.e, internal rays ${\cal R}$ satisfying $P^{\circ n}({\cal R})={\cal R}$. 
They correspond in the B\"ottcher coordinate to the segments 
$\{re^{{2\pi ki} \over {{\cal D}-1}}:r \in [0,1),\hbox{ } k=0,\dots,{\cal D}-2\}$. 
\smallskip
What is important to note here, is that the same construction is valid for all elements in the critical cycle. 
Note that if we choose a coordinate $\phi_\omega$ in which the internal ray ${\cal R}$
corresponds to the real segment $[0,1)$, 
then we can choose in a unique way a coordinate $\phi_{P(\omega)}$ (at $P(\omega)$) for which $P({\cal R})$ corresponds to $[0,1)$. 
Furthermore in this case 
$$
\phi_{P(\omega)}(P(z))=(\phi_\omega(z))^{deg_\omega P},
$$
where $deg_\omega P$ is the local degree of $P$ at $\omega$ 
(for more details see [DH1, Chapter 4, Proposition 2.2]).
\bigskip 
{\bf 1.2 Lemma.} 
{\it If a critical point $z$ belongs to a critical cycle of period $n=n_z$, 
then $P^{\circ n}|_{\overline{U(z)}}$ (which has degree say ${\cal D}_z>1$) 
has 
exactly ${\cal D}_z-1$ different fixed points 
in the boundary $\partial U$ of this component $U(z)$ respect to this return map. 
Furthermore, all external rays that land at such points have period exactly $n$.} 
\medskip
{\bf Proof.} 
The first part is well known. 
For the second, 
we consider near this periodic point segments of all the external rays which land there, 
together with the internal ray joining this point to the center $z$. 
The cyclic order of these analytic arcs must be preserved under iteration. 
The result thus follows easily. 
\endofproof
\bigskip
{\bf 1.3 Supporting arguments.} 
Given a Fatou component $U$ and a point $p \in \partial U$, 
there are only a finite number of external rays $R_{\theta_1},\dots,R_{\theta_k}$ landing at $p$. 
These rays divide the plane in $k$ regions.
We order the arguments of these rays in counterclockwise cyclic order 
$\{\theta_1,\dots,\theta_k\}$, 
so that $U$ belongs to the region determined by 
$R_{\theta_1}$ and $R_{\theta_2}$  
($\theta_1=\theta_2$ if there is a single ray landing at $p$). 
The argument $\theta_1$ (respectively the ray $R_{\theta_1}$) is by definition the 
{\it (left) supporting argument (respectively the (left) supporting ray) of the Fatou component U}. 
In a completely analogous way we can define right supporting rays. 
Note that an argument supports at most one Fatou component 
(compare [DH1, Chapter VII.4]).
Furthermore, by definition, 
given a Fatou component $U$, 
for every boundary point $p$ there is an external ray landing at $p$, 
and therefore a supporting ray for $U$. 
\medskip
{\bf 1.4 Extended Rays.} 
Given an external ray $R_\theta$ supporting the Fatou component $U(z)$ with center $z$, 
we extend $R_\theta$ by joining its landing point with $z$ by an internal ray, 
and call this set an {\it extended ray $\hat R_\theta$ with\break 
argument $\theta$}. 
\bigskip 
{\bf 1.5 Example.} 
Consider the postcritically finite polynomial $P_c(z)=z^2+c$ 
(where the critical value $c \approx -0.12256117+0.74486177i$ satisfies $c^3+2c^2+c+1=0$). 
The rays with argument $1/7,2/7,4/7$ all land at the same fixed point. 
But $R_{4/7}$ is the only ray landing at this point, 
which supports the critical component. 
(Compare Figure 1.1.)
\bigskip
\centerline{
\insertRaster 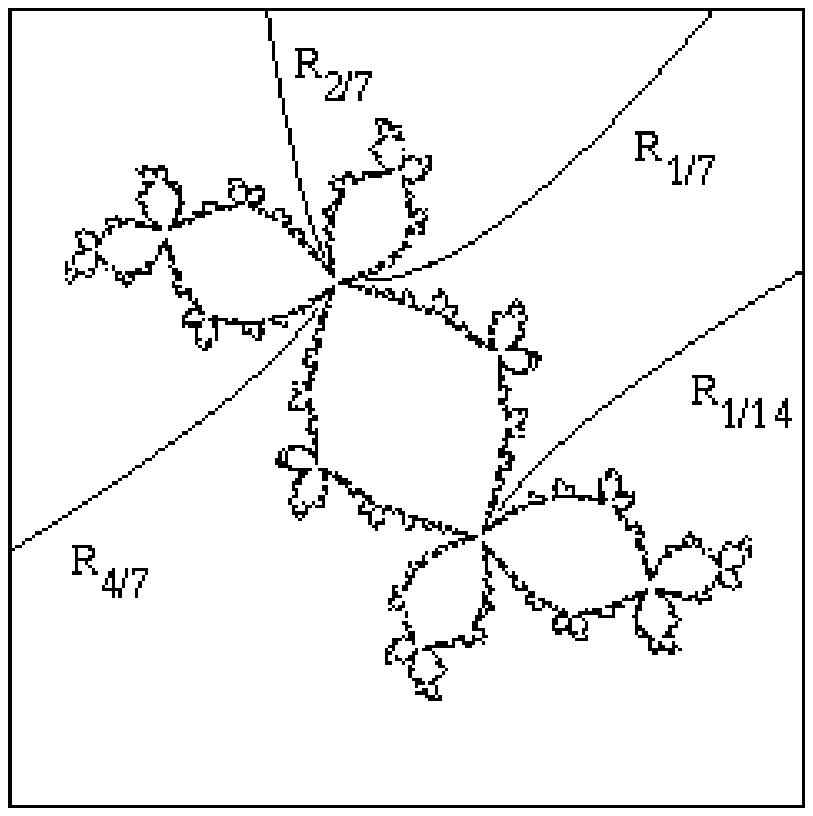 pixels 256 by 240 scaled 675
}
\medskip
\centerline{\it Figure 1.1}
\bigskip
\insec
\centerline{\medfont 2. Construction of Critically Marked Polynomials.} 
\insec 
Given a postcritically finite polynomial $P$, 
we associate to every critical point a finite subset of ${\bf Q/Z}$ and construct a 
{\it critically marked polynomial} 
$(P, {\bf \cal F} = \{{\cal F}_1, \dots, {\cal F}_{n_F}\}$, 
${\bf \cal J}=\{{\cal J}_1,\dots,{\cal J}_{n_J}\})$. 
Here ${\cal F}_k$ would be the set of arguments associated with the critical point $z^F_k$ in the Fatou set, 
and ${\cal J}_k$ would be the set associated with the critical point $z^J_k$ in the Julia set. 
The number of elements in these finite sets would be equal to the local degree of the associated critical points. 
We remark that given a polynomial its critical marking is not necessarily unique. 
Also note that one of these two families will be empty if there are no critical points in the Fatou or Julia set.
In the following definition we will always work with left supporting rays. 
We remark that we could equally well work with the right analogue, 
but there must be the same choice throughout. 
Also, multiplication by $d$ modulo $1$ in {\bf R/Z} will be denoted by $m_d$. 
\bigskip
{\bf 2.1 Construction of ${\cal F}_k$.} 
First we consider the case in which a given Fatou critical point $z=z^F_k$ is periodic. 
Let $z=z^F_k \mapsto P(z) \mapsto ...\break 
\mapsto P^{\circ n}(z)=z$ be a critical cycle of period $n$ and degree ${\cal D}_z>1$ 
(compare $\S$1.1). 
We construct the associated set ${\cal F}_\ell$ for every critical point $z^{\cal F}_\ell$ in the cycle simultaneously. 
Denote by $d_z$ the local degree of $P$ at $z$. 
We pick any periodic point $p_z \in \partial U(z)$ of period dividing $n$ 
(which is not critical because it is periodic and belongs to the Julia set $J(P)$) 
and consider the supporting ray $R_\theta$ for this component $U(z)$ at $p_z$. 
Note that this choice naturally determines a periodic supporting ray for every Fatou component in the cycle. 
The period of this ray is exactly $n$ (compare Lemma 1.2). 
Given this periodic supporting ray $R_\theta$, 
we consider the $d_z$ supporting rays for this same component $U(z)$ that are inverse images of $P(R_\theta)=R_{m_d(\theta)}$. 
The set of arguments of these rays is defined to be ${\cal F}_k$. 
Keeping in mind that a preferred periodic supporting ray has been already chosen, 
we repeat the same construction for all critical points in this cycle. 
Note that as the cycle has critical degree ${\cal D}_z$, 
we can produce ${\cal D}_z-1$ different possible choices for ${\cal F}_k$. 
If ${\cal F}_k$ is the set associated with the periodic critical point $z_k$, 
there is only one periodic argument in ${\cal F}_k$ 
(namely $\theta$ as above), 
we call this angle {\it the preferred supporting argument associated with $z^F_k$.} 
Note that by definition, the period of $z_k^F$ 
equals the period of the associated preferred periodic argument. 
\smallskip
Otherwise, if $z=z^F_k$ of degree $d_z>1$, 
is a non periodic critical point in the Fatou set $F(P)$, 
there exists a minimal $n>0$ for which $w=P^{\circ n}(z)$ is critical. 
If $w$ has associated a preferred supporting ray $R_\theta$ 
(at the beginning only periodic critical points do), 
then $P^{-n}(R_\theta)$ contains exactly $d_z$ rays which support this Fatou component $U(z)$. 
The set of arguments of these rays is defined to be ${\cal F}_k$. 
We pick any of those and call it the {\it preferred supporting argument associated with $z$}. 
We continue this process for all Fatou critical points.
\medskip
{\bf 2.2 Construction of ${\cal J}_k$.} 
Given $z=z_k^J$ (a critical point in $J(P)$) of degree $d_k>1$, 
we distinguish two cases. 
If the forward orbit of $z$ contains no other critical point, 
we have that for some $\theta$ (usually non unique) $R_\theta$ lands at $P(z)$. 
Now $P^{-1}(R_\theta)$ consists of $d$ different rays, 
among them exactly $d_k$ land at $z$. 
Define ${\cal J}_k$ as the set of arguments of these rays, 
and choose a {\it preferred ray}. 
Otherwise, $z$ will map in $n \ge 1$ iterations to a critical point, 
which we assume to have associated a preferred ray $R_{\theta}$. 
In the $n^{th}$ inverse image $P^{-n}(R_\theta)$ of this preferred ray, 
there are $d_k$ rays which land at $z$. 
The set of arguments of these rays is defined to be ${\cal J}_k$. 
Again we pick one of those to be preferred, 
and continue until every critical point has an associated set. 
\bigskip 
The critical marking itself gives information about 
how many iterates are needed for a given critical point to become periodic. 
For example we have the following lemma. 
\smallskip
{\bf 2.3 Lemma.} 
{\it Let $\gamma$ be a preferred supporting argument in the set ${\cal F}_k$\break 
(respectively in ${\cal J}_k$). 
Then the multiple 
$m^{\circ n}_d(\gamma)$ (with $n \ge 1)$ is periodic but $m_d^{\circ n-1}(\gamma)$ is not 
if and only if 
$z^F_k$ (respectively $z^J_k$) 
falls in exactly $n$ iterations into a periodic orbit.} 
\medskip
{\bf Proof.} 
This clearly follows from the construction. 
\endofproof 
\medskip
The importance of the above construction is stated in the following theorem. 
The proof will be given in Chapter III (compare also Theorem 3.9). 
\medskip
{\bf 2.4 Theorem.}
{\it Every centered monic postcritically finite polynomial P has a critical marking $(P,{\cal F},{\cal J})$. 
This marking determines the polynomial P in the following sense: 
if $(P,{\cal F},{\cal J})$ and $(Q,{\cal F},{\cal J})$ are critically marked polynomials, then $P=Q$. 
In other words, 
two monic centered post-critically finite polynomials with the same critical marking $({\cal F},{\cal J})$ must be equal.} 
\medskip
{\bf Remark.} 
Note that the construction of associated sets was done in several steps. 
We first complete the choice for all critical cycles, 
and then proceed backwards. 
In both the Fatou and Julia set cases we will have to make decisions at several stages of the construction. 
Such decisions will affect the choice of the marking for all critical points found in the backward orbit of these starting ones. 
Each time that this kind of construction is made, 
we will informally say that it is a {\it hierarchic selection}. 
We encourage the reader to take a look at the examples in Section 4. 
\vfill\eject
\insec
\centerline{\medfont 3. The Combinatorics of Critically Marked Polynomials.} 
\insec 
In order to analyze which properties the families $({\bf \cal F},{\bf \cal J})$ satisfy, 
it is convenient to introduce some combinatorial notation. 
\bigskip
{\bf 3.1 Definitions.} 
We say that a subset ${\Lambda} \subset {\bf R/Z}$ is a {\it (d-)preargument set} 
if $m_d({\Lambda})$ is a singleton. 
For technical reasons we will always assume that $\Lambda$ contains at least two elements. 
If all elements of $\Lambda$ are rational, we say that $\Lambda$ is a {\it rational preargument set}. 
It follows by construction that whenever $(P,{\bf {\cal F},{\cal J}})$ is a marked polynomial, 
all the sets ${\cal J}_k$, and ${\cal F}_l$ are rational $d$-preargument sets.
\smallskip
Consider now a family ${\bf \Lambda}=\{\Lambda_1,\dots,\Lambda_n\}$ of finite subsets of the circle ${\bf R/Z}$. 
The family ${\bf \Lambda}$ determines the {\it family union} set ${\bf \Lambda}^\cup =\bigcup \Lambda_i$. 
We say that any $\lambda \in {\bf \Lambda}^\cup$ is an {\it element of the family ${\bf \Lambda}$}. 
Furthermore, we can say that it is a periodic or preperiodic element of the family if it is so with respect to $m_d$. 
The set of all periodic elements in the family union will be denoted by 
${\bf \Lambda^\cup_{per}}$.
\bigskip
{\bf 3.2 Hierarchic Families.} 
We say that a family ${\bf \Lambda}$ is {\it hierarchic} if 
for any elements in the family $\lambda,\lambda' \in {\bf \Lambda}^\cup$, 
whenever 
$m_d^{\circ i}(\lambda),m_d^{\circ j}(\lambda') \in \Lambda_k$ for some $i,j > 0$ 
then 
$m_d^{\circ i}(\lambda)=m_d^{\circ j}(\lambda')$.
(This is useful if we think of a dynamically preferred element in each $\Lambda_k$).
\bigskip
{\bf 3.3 Linkage Relations.} 
We will say that two subsets $T$ and $T'$ of the circle $\bf R/Z$ are unlinked if 
they are contained in disjoint connected subsets of $\bf R/Z$, 
or equivalently, 
if $T'$ is contained in just one connected component of the complement ${\bf R/Z}-T$. 
(In particular $T$ and $T'$ must be disjoint.) 
If we identify $\bf R/Z$ with the boundary of the unit disk, 
an equivalent condition would be that the convex closures of these sets are pairwise disjoint. 
If $T$ and $T'$ are not unlinked then either $T \cap T'\ne \emptyset$ or there are elements $\theta_1,\theta_2 \in T$ and 
$\theta'_1,\theta'_2 \in T'$ such that the cyclic order can be written
$\theta_1, \theta'_1, \theta_2, \theta'_2,\theta_1$. 
In this second case we say that $T$ and $T'$ are {\it linked}.
More generally a family ${\bf \Lambda}=\{\Lambda_1,\dots,\Lambda_n\}$ is an {\it unlinked family} if 
$\Lambda_1,\dots,\Lambda_n$ are pairwise unlinked.
Alternatively, each $\Lambda_i$ is completely contained in a component of 
${\bf R/Z}-\Lambda_j$ for all $j \ne i$.
\medskip
The preceding definition has its motivation in the description of the dynamics of external rays for a polynomial map. 
Suppose the external rays $R_{\theta_i},R_{\psi_i}$ land at $z_i$ for $i=1,2$. 
If $z_1 \ne z_2$ then the sets $\{\theta_1,\psi_1\},\{\theta_2,\psi_2\}$ are unlinked, 
for otherwise the rays will cross each other.
The same argument applies if we consider rays supporting Fatou components. 
But if we analyze linkage relations arising from rays supporting a Fatou component and rays that land at some point, 
we may get minor problems. 
Anyway, it is easy to see that even in this case the associated sets of arguments will be `almost' unlinked. 
(Compare condition (c.2) and as well as Proposition 3.8 below.) 
\bigskip 
{\bf 3.4 Weak linkage relations.} 
Consider two families ${\bf {\cal F}}=\{{\cal F}_1,\dots,{\cal F}_n\}$ and 
${\bf {\cal J}}=\{{\cal J}_1,\dots,{\cal J}_m\}$; 
we say that ${\bf \cal J}$ is {\it weakly unlinked to ${\bf \cal F}$ in the right} 
if we can chose arbitrarily small 
$\epsilon>0$ so that the family 
$\{{\cal F}_1,\dots,{\cal F}_n,\break 
{\cal J}_1-\epsilon,\dots,{\cal J}_m-\epsilon\}$ is unlinked. 
(Here 
${\Lambda}-\epsilon=\{\lambda-\epsilon \hbox{ ({\it mod 1}) : }\lambda \in {\Lambda}\}$.)
In particular each family should be unlinked. 
Note that the definition allows empty families. 
To simplify notation we will simply say that 
{\it ``${\cal F}$ and ${\cal J}^-$ are unlinked".} 
\bigskip
{\bf 3.5 Formal Critical Portraits.}
Consider families ${\bf {\cal F}}=\{{\cal F}_1,\dots,{\cal F}_n\}$ and 
${\bf {\cal J}}=\{{\cal J}_1,\dots,{\cal J}_m\}$ of rational ({\it d-})prearguments. 
We say that the pair $({\bf {\cal F},{\cal J}})$ is a {\it degree d formal critical portrait} if the following conditions are satisfied. 
\smallskip
{\it 
(c.1) $d-1=\sum (\#({\bf \cal F}_k)-1)+\sum (\#({\bf \cal J}_l)-1)$
\tskip
(c.2) ${\cal F}$ and ${\cal J}^-$ are unlinked. 
\tskip
(c.3) Each family is hierarchic.
\tskip
(c.4) Given $\gamma \in {\bf {\cal F}^\cup}$, there is an $i>0$ such that $m_d^{\circ i}(\gamma) \in {\bf {\cal F}^\cup_{per}}$.
\tskip
(c.5) No $\theta \in {\bf {\cal J}^\cup}$ is periodic.}
\medskip
This set of conditions represent the simplest conditions satisfied 
by the critical marking of a postcritically finite polynomial. 
Condition (c.1) says that we have chosen the right number of arguments. 
Condition (c.2) means that the rays and extended rays determine sectors which do not cross each other, 
and that ${\bf \cal F}$ was constructed from arguments of left supporting rays. 
This reflects our decision to chose the supporting arguments as the rightmost possible argument of an external ray. 
Condition (c.3) reflects our choice of dynamically preferred rays. 
Condition (c.4) indicates that arguments in ${\bf \cal F}$ are related to Fatou critical points. 
Condition (c.5) indicates that arguments in ${\bf \cal J}$ are related to Julia set critical points.
Unfortunately there are formal critical portraits which do not correspond to a postcritically finite polynomial 
(compare Example II.2.8). 
In order to state necessary and sufficient conditions we need to study the 
dynamically defined partitions of the unit circle determined by these elements. 
\bigskip
{\bf 3.6.} 
Given two families ${\bf {\cal F},{\cal J}}$ as above, 
we form a partition ${\bf \cal P}=\{L_1,\dots,L_d\}$ 
of the unit circle minus a finite number of points 
${\bf R/Z}-{\bf {\cal F}^\cup} - {\bf {\cal J}^\cup}$, in the following way. 
We consider two points $t,t' \in {\bf R/Z}-{\bf {\cal F}^\cup}-{\bf {\cal J}^\cup}$. 
By definition, $t,t'$ are {\it unlink equivalent} 
if they belong to the same connected component of ${\bf R/Z}- {\cal F}_i$ and ${\bf R/Z}- {\cal J}_j$, for all possible $i,j$.
Let $L_1,\dots,L_d$ be the resulting unlink equivalence classes with union ${\bf R/Z}-{\bf {\cal F}^\cup} - {\bf {\cal J}^\cup}$. 
It is easy to check that each $L_p$ is a finite union of open intervals with total length $1/d$.
\smallskip
Each element $L_i \in {\cal P}$ of the partition 
is a finite union $L_i=\cup (x_j,y_j)$ of open connected intervals.  
We also define the sets $L_i^+=\cup [x_j,y_j)$ and\break 
$L_i^-=\cup (x_j,y_j]$. 
It is easy to see that both 
${\bf \cal P}^+=\{L_1^+,\dots,L_d^+\}$ and\break 
${\bf \cal P}^-=\{L_1^-,\dots,L_d^-\}$ 
are partitions of the unit circle. 
As every $\theta \in {\bf R/Z}$ belongs to exactly one set $L_k^+$, 
we define its {\it right address} $A^+(\theta)=L_k$. 
In an analogous way we define the {\it left address $A^-(\theta)$ of $\theta$}. 
We associate to every argument $\theta \in {\bf R/Z}$  
a {\it right symbol sequence} $S^+(\theta)=(A^+(\theta),A^+(m_d(\theta)),\dots)$ and 
a {\it left symbol sequence} 
$S^-(\theta)=(A^-(\theta),A^-(m_d(\theta)),\dots)$. 
Note that for all but a countable number of arguments $\theta \in {\bf R/Z}$ 
(namely the arguments present in the families and their iterated inverses), 
the left $S^-(\theta)$ and the right $S^+(\theta)$ symbol sequences coincide. 
By $S(\theta)$ will be meant either (left or right) symbol sequence. 
\bigskip
{\bf 3.7 Admissible Critical Portraits.}
Let ${\bf \cal F}=\{{\cal F}_1,\dots,{\cal F}_n\}$ and\break 
${\bf \cal J}=\{{\cal J}_1,\dots,{\cal J}_m\}$ 
be two families of rational ({\it d-})prearguments. 
We say that $({\bf {\cal F},{\cal J}})$ is a 
{\it degree d admissible critical portrait} if 
$({\bf {\cal F},{\cal J}})$ is a degree $d$ formal critical portrait and 
the following two extra conditions are satisfied. 
\smallskip
{\it 
(c.6) Let $\gamma \in {\bf {\cal F}^\cup_{per}}$ and $\lambda \in {\bf R/Z}$, 
then $\lambda=\gamma$ if and only if\break 
$S^+(\gamma)=S^+(\lambda)$.
\tskip 
(c.7) Let $\theta \in {\cal J}_l$ and $\theta' \in {\cal J}_k$. 
If for some i, $S^-(m_d^{\circ i}(\theta))=S^-(\theta')$, then $m_d^{\circ i}(\theta) \in {\cal J}_k$.}
\bigskip
{\bf 3.8 Proposition.}
{\it If $(P,{\cal F},{\cal J})$ is a critically marked polynomial, then $({\cal F},{\cal J})$ is an admissible critical portrait.} 
\medskip
Condition (c.6) indicates that arguments in ${\cal F}_l$ must support Fatou components. 
Condition (c.7) indicates that different elements in the family ${\cal J}$ are associated with different critical points. 
The proof of this proposition will be given in Section II.2. 
\bigskip
Now we can state the main result for critically marked polynomials as follows 
(the proof of this theorem will be given in Chapter III). 
\medskip 
{\bf 3.9 Theorem.}
{\it Let $({\bf {\cal F},{\cal J}})$ be a degree d admissible critical portrait. 
Then there is a unique monic centered postcritically finite polynomial P, 
with critical marking $(P,{\bf {\cal F},{\cal J}})$.}
\bigskip 
Now we should ask if conditions {\it (c.1)-(c.7)} represent a finite amount of information to be checked. 
This question is answered in a positive way by the following lemma. 
The proof would be given in Section II.1. 
\medskip
{\bf 3.10 Lemma.}
{\it Suppose $\theta$ and $\theta'$ have the same periodic left (or right) symbol sequence. 
Then $\theta$ and $\theta'$ are both periodic and of the same period.} 
\bigskip
In particular condition $(c.6)$ can be replaced by condition $(c.6)'$: 
\medskip
{\it 
(c.6)$'$ Let $\gamma \in {\bf {\cal F}^\cup_{per}}$ and let $\lambda$ have the same period as $\gamma$, 
then $\lambda=\gamma$ if and only if $S^+(\gamma)=S^+(\lambda)$. 
} 
\bigskip
{\bf 3.11.} The next question that we ask is what kind of information 
about the Julia set can be gained by looking carefully into the combinatorics. 
For example, 
if can we determine if two rays land at the same point by only looking at their arguments. 
In fact, left symbol sequences 
convey all the information necessary to effectively decide whether two rays land at the same point or not. 
This is done as follows. 
Suppose ${\cal J}_i=\{\theta_1,...,\theta_k\} \in {\cal J}$ with corresponding left symbol sequences 
$S^-(\theta_1),\dots,S^-(\theta_k)$. 
As we expect the rays with those arguments to land at the same critical point, 
we declare them ($i$-)equivalent; 
i.e, we write $S^-(\theta_\alpha) \equiv_i S^-(\theta_\beta)$. 
Then we set $\theta \approx \theta'$ either if 
$S^-(\theta)=S^-(\theta')$ or 
there is an $n \ge 0$ such that $A^-(m_d^{\circ j}(\theta))=A^-(m_d^{\circ j}(\theta'))$ for all $j < n$ and 
$S^-(m_d^{\circ n}(\theta)) \equiv_i S^-(m_d^{\circ n}(\theta'))$ for some $i$. 
This relation $\approx$ is not necessarily an equivalence relation, because transitivity may fail. 
To make this into an equivalence relation we say that $\theta \sim_l \theta'$ 
if and only if 
there are arguments $\lambda_0=\theta,\lambda_1,\dots,\lambda_m=\theta'$, 
such that $\lambda_0 \approx \dots \approx \lambda_m$. 
The importance of this equivalence relation is shown by the following proposition. 
The proof will be given in Chapter II. 
\bigskip
{\bf 3.12 Proposition.} 
{\it Let $(P,{\bf {\cal F},{\cal J}})$ be a critically marked polynomial. 
Then $R_\theta$ and $R_{\theta'}$ land at the same point if and only if $\theta \sim_l \theta'$.} 
\bigskip
{\bf 3.13 Corollary.} 
{\it 
The symbol sequence $S^-(\theta)$ is a periodic sequence of period $m$ 
if and only if 
the landing point of the ray $R_\theta$ has period $m$.} 
\endofproof
\bigskip
We proceed now to give a very brief description of Chapters II and III 
which are devoted to critical portraits. 
In Chapter II we will work in more detail the combinatorics of critical portraits. 
We will also translate several of our results to the corresponding Julia set. 
In Chapter III we give the proof of the Realization Theorem for Critical Portraits. 
Appendix A deals with the relevant part for our use of 
Thurston's theory of postcritically finite rational maps. 
We state Thurston's Theorem in a more general form. 
Namely, we include the possibility of additional periodic or preperiodic orbits. 
The proof given in [DH2] extends to this formulation. 
\bigskip
We now give a brief comparison between our work and that by 
Bielefeld, Fisher and Hubbard. 
Of course there is a big overlap in both expositions, 
and we have, when possible, referred to the original proofs. 
In Chapter II, we analyze the combinatorics of the marking and the proofs 
follow the same lines as those in [BFH]. 
Chapter III is essentially different and we have stated without proofs those results in [BFH] 
which apply to our case. 
The main point here is that as Levy cycles can not involve any 
`preperiodic element' of the topological polynomial, 
new preperiodic arguments should be introduced artificially. 
We call them {\it special arguments}. 
Finally, we still have to prove that the recovered polynomial admits the specified critical marking. 
Our method of proof is more delicate because new difficulties are involved. 
\bigskip
\bigskip
{\bf Acknowledgement.} 
We will like to thank John Milnor for helpful discussions and suggestions. 
Some of the arguments are in its final formulation thanks to him. 
We will also want to thank (among others) to 
Ben Bielefeld and John Hubbard for discussions at different stages of the preparation of this work. 
Most of the figures were constructed using a program of Milnor. 
Also, we want to thank the Geometry Center, University of Minnesota
and Universidad Cat\'olica del Per\'u 
for their material support. 
\bigskip

\insec
\centerline{\medfont 4. Examples.}
\insec
{ 
We will illustrate with examples the definitions of the previous sections. 
We will try to isolate and illustrate all possible complications. 
Of course, the worst possible examples will involve several of these at the same time.}  
\bigskip
{\bf 4.1 The rabbit.} 
(See Figure 1.1.) 
Once again consider the degree two polynomial $P_c(z)=z^2+c$ with $c \approx -0.12256117+0.74486177i$. 
The Fatou critical point $z=0$ has a period 3 orbit under iteration. 
Therefore $P^{\circ 3}$ restricted to the critical component is a degree 2 cover of itself. 
It follows that the map $P^{\circ 3}$ has a unique fixed point in the boundary of this critical Fatou component. 
As noticed above, among the three rays $R_{1/7},R_{2/7},R_{4/7}$ landing at this fixed point, 
only the ray $R_{4/7}$ supports the critical component. 
By the definition of marking, 
we must look for the other ray that supports this component and maps to $P(R_{4/7})=R_{1/7}$. 
This ray can only be $R_{1/14}$. 
Thus, we have constructed a marking for $P$. 
In this case ${\bf \cal F}=\{{\cal F}_1\}$ and ${\bf \cal J}=\emptyset$, 
where ${\cal F}_1=\{4/7,1/14\}$. 
\smallskip
It is important to note that we were looking for a fixed point of $P^{\circ 3}$ restricted to the boundary of the critical Fatou component. 
Such a fixed point for $P^{\circ 3}$ turned out to be a fixed point for $P$ as well, 
but the rays landing there have period equal 3.
\bigskip 
{\bf 4.2 The Ulam-von Neumann map.} 
We consider now the strictly preperiodic case. 
Let $P(z)=z^2-2$, and note that the orbit of the critical point $z=0$ is 
$0 \mapsto -2 \mapsto 2 \mapsto 2 \dots$. 
Only the external ray $R_0$ lands at $z=2$, 
and therefore only the ray $R_{1/2}$ lands at $z=-2$. 
Both $R_{1/4},R_{3/4}$ land at the critical point $z=0$, 
and map to $R_{1/2}$ under $P$. 
In this case the marking is ${\bf \cal F}=\emptyset$ and ${\bf \cal J}=\{{\cal J}_1\}$, 
where ${\cal J}_1=\{1/4,3/4\}$. 
\bigskip
{\bf 4.3 Preperiodic case: two possible choices.} 
(See Figure 1.2.) 
Consider the degree two polynomial $P_c(z)=z^2+c$ where $c \approx-1.5436891$ is the only negative solution of the equation $c^3+2c^2+2c+2=0$. 
\smallskip
In this case the critical point $z=0$ has orbit $0 \mapsto c \mapsto c^2+c\break 
\mapsto -(c^2+c) \mapsto -(c^2+c)$. 
The rays $R_{1/3},R_{2/3}$ both land at the fixed point $z=-(c^2+c)$, and are interchanged by $P_c$. 
At $z=c^2+c$, the rays $R_{1/6},R_{5/6}$ land. 
At the critical value $z=c$, $R_{5/12},R_{7/12}$. 
In this way, we can get two different markings ${\bf \cal F}=\emptyset$, ${\bf \cal J}=\{{\cal J}_1\}$, 
where ${\cal J}_1=\{5/24,17/24\}$ corresponds to the choice of the {\it critical ray} $R_{5/12}$, 
and  ${\cal J}_1=\{7/24,19/24\}$ to the choice of $R_{7/12}$. 
\smallskip
In either case we can read from the marking that the critical point needs three iterations to become periodic. 
The exact period however can not be read immediately from this data. 
(Compare Corollary 3.13.) 
\bigskip
\centerline{
\insertRaster 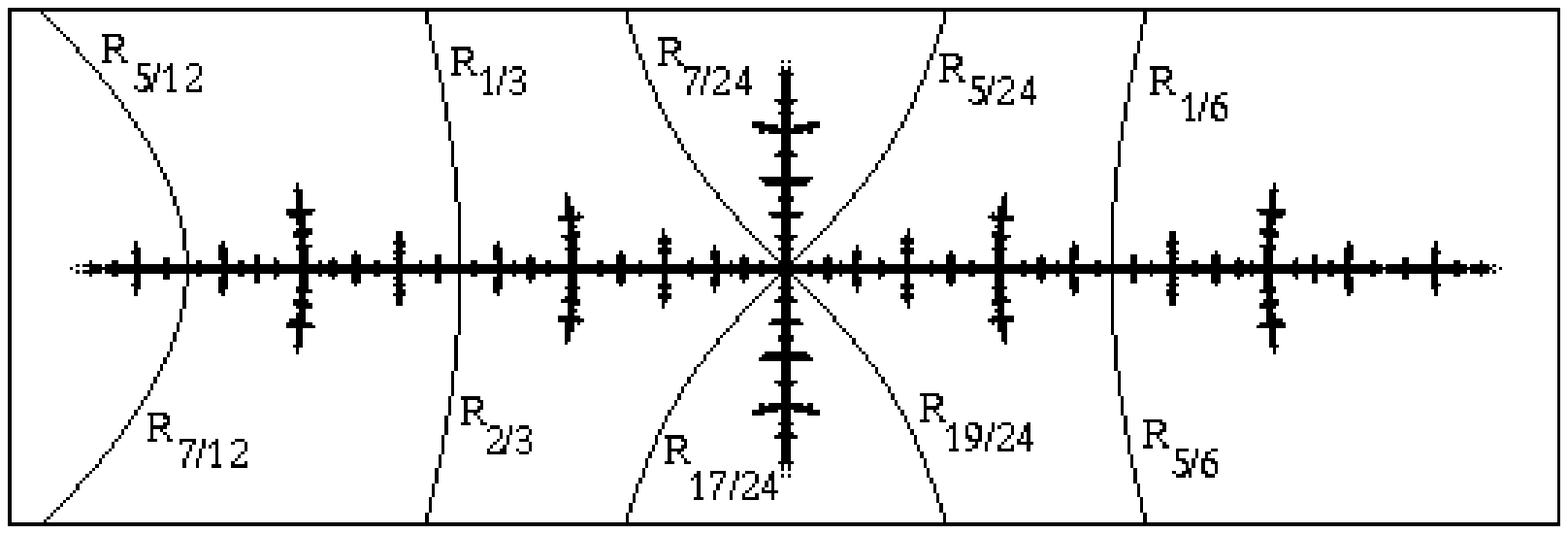 pixels 480 by 160 scaled 725
}
\medskip
\centerline{\it Figure 1.2}
\bigskip
{\bf 4.4 Non trivial critical cycle.} 
(See Figure 1.3.) 
Consider the degree 3 polynomial $P(z)=z^3-{3 \over 2}z$. 
The critical points satisfy $z^2=1/2$, and is easy to see that they are interchanged by $P$ 
(i.e, if $a$ is a critical point then $P(a)=-a$). 
In each of the critical Fatou components the map $P^{\circ 2}$ is a degree 4 
(the product of the degrees of the cycle!) 
covering of itself. 
In this way, there must be in the boundary of each component 3 ($=4-1$) possible choices of periodic points. 
One of those fixed points ($z=0$) belongs to the boundary of both components. 
The rays landing at $z=0$ are $R_{1/4}$ and $R_{3/4}$, 
and each one supports exactly one of the Fatou critical components. 
The period 2 rays that support the {\it `rightmost'} component are $R_{3/4},R_{7/8},R_{1/8}$ 
(their respective images $R_{1/4},R_{5/8},R_{3/8}$ support the other). 
Therefore, the choice of a periodic supporting ray for one component, 
forces the choice of its image for the other. 
\bigskip
\medskip
\centerline{
\insertRaster 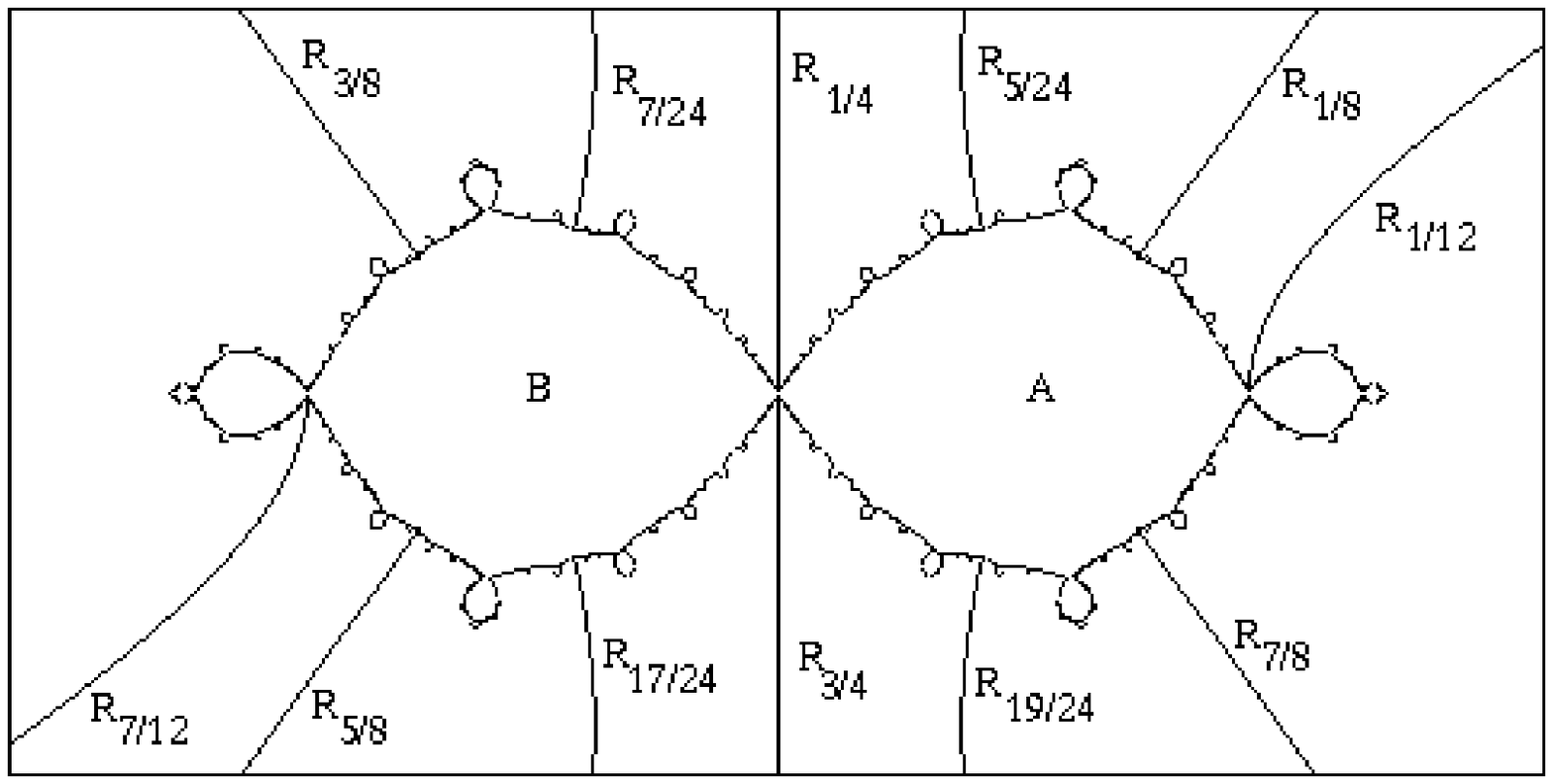 pixels 480 by 240 scaled 550
}
\medskip
\centerline{\it Figure 1.3}
\bigskip 

This polynomial has exactly three markings, all of type 
${\bf \cal F}=\{{\cal F}_A,{\cal F}_B\}$,\break 
${\bf \cal J}=\emptyset$. 
The periodic supporting rays are listed on the left. 
\bigskip
\hskip 0.05in {\it Component A} \hskip 0.35in {\it Component B} \hskip 0.3in
${\cal F}_A$ \hskip 0.9in ${\cal F}_B$ 
\tskip
\hskip 0.3in $R_{3/4}$ \hskip 1in $R_{1/4}$ \hskip 0.3in $\{3/4, 1/12\}$
\hskip 0.4in $\{1/4, 7/12\}$
\tskip
\hskip 0.3in $R_{7/8}$ \hskip 1in $R_{5/8}$ \hskip 0.3in $\{7/8, 5/24\}$
\hskip 0.4in $\{5/8, 7/24\}$
\tskip
\hskip 0.3in $R_{1/8}$ \hskip 1in $R_{3/8}$ \hskip 0.3in $\{1/8,19/24\}$
\hskip 0.35in $\{3/8,17/24\}$
\bigskip
The question is now, 
why can we not take ${\cal F}_A=\{3/4, 1/12\}$ and\break 
${\cal F}_B=\{3/8,17/24\}$ as a marking? 
This is forbidden by the rules of $\S 3$ 
since $3/4$ and $3/8$ do not belong to the same cycle. 
A good reason for this rule is given in the next example. 
\bigskip
\centerline{
\insertRaster 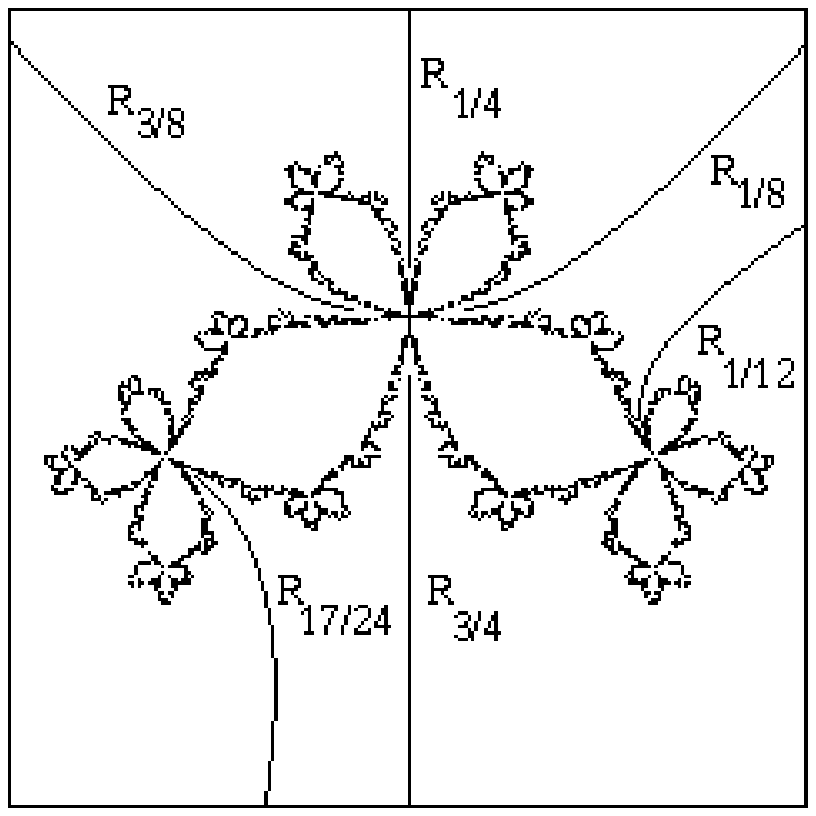 pixels 256 by 240 scaled 675
}
\medskip
\centerline{\it Figure 1.4}
\bigskip
{\bf 4.5 Bad choice, wrong polynomial.} 
(See Figure 1.4.)
There is a polynomial with marking 
${\bf \cal F}=\{{\cal F}_A,{\cal F}_B\},{\bf \cal J}=\emptyset$, 
where ${\cal F}_A=\{3/4, 1/12\}$, ${\cal F}_B=\{3/8,17/24\}$. 
But it is not the one in Example 4.4. 
\smallskip
For the polynomial $P(z)=z^3+az+b$ (where $a=-0.75, b \approx 0.661438i$), 
the rays $R_{3/4},R_{1/8},R_{1/4},R_{3/8}$, 
land at a fixed point which belongs to the boundary of the four periodic Fatou components. 
Those components are associated pairwise in cycles, so we have two disjoint degree 2 cycles. 
Only $R_{3/4}$ and $R_{3/8}$ support critical components. 
It follows easily that this polynomial has a unique marking. 
\bigskip
{\bf 4.6 Hierarchic choice.} 
(See Figure 1.5.) 
Consider now the polynomial 
$P(z)={\sqrt{2}}(z^2-1)^2$, with critical points $z=0,\pm1$. 
The orbit of the critical points is 
$\pm1 \mapsto 0 \mapsto \sqrt{2} \mapsto \sqrt{2}$. 
At the fixed point $z=\sqrt{2}$ only the ray $R_0$ lands. 
At $z=0$, $R_{1/4}$ and $R_{3/4}$.
At $z=1$, $R_{1/16}, R_{3/16}, R_{13/16}, R_{15/16}$. 
At $z=-1$, $R_{5/16}, R_{7/16}, R_{9/16}, R_{11/16}$. 
In this case the marking will not be unique and will depend in the choice of the preferred ray at $z=0$. 
\bigskip
\centerline{
\insertRaster 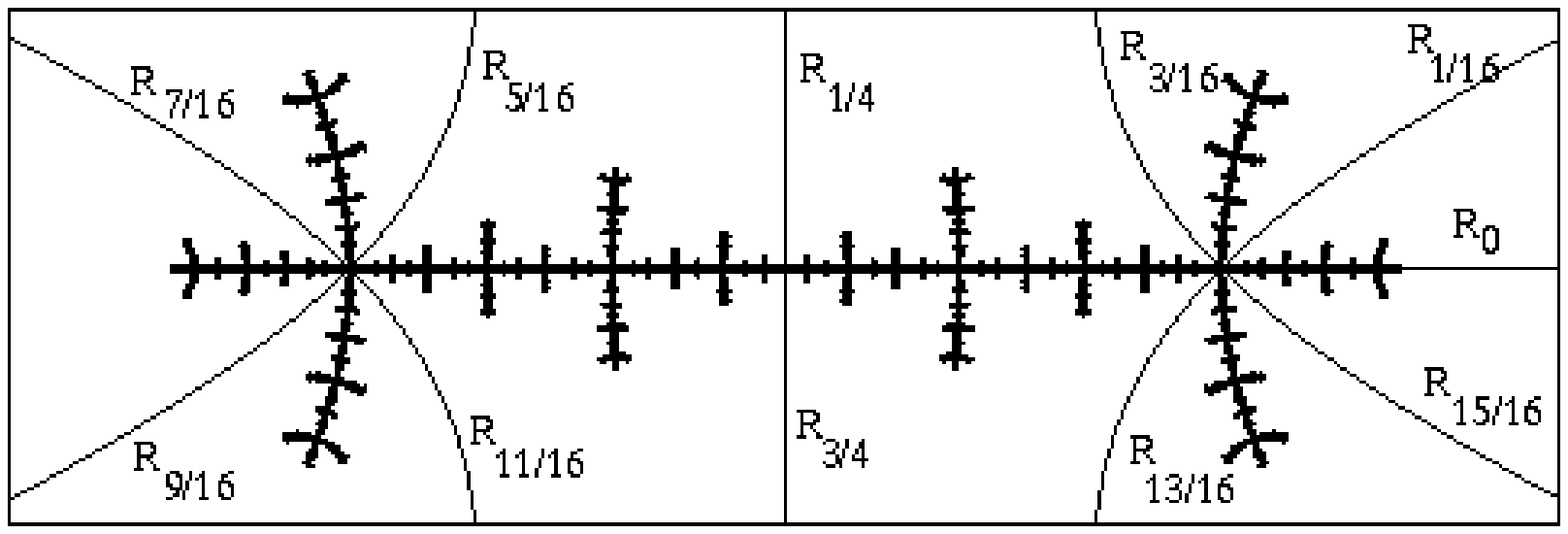 pixels 494 by 160 scaled 725
}
\medskip
\centerline{\it Figure 1.5}
\bigskip
The marking will be of the form 
${\bf \cal F}=\emptyset,{\bf \cal J}=\{{\cal J}_{z=0},{\cal J}_{z=1},{\cal J}_{z=-1}\}$. 
\bigskip 
\hskip 0.15in ${\cal J}_{z=0}$ \hskip 0.5in {\it preferred ray} \hskip 0.6in 
${\cal J}_{z=1}$ \hskip 0.7in ${\cal J}_{z=-1}$
\tskip
\hskip 1.05in {\it at $z=0$}
\tskip
$\{1/4,3/4\}$ \hskip 0.55in $R_{1/4}$ \hskip 0.6in $\{1/16,13/16\}$ \hskip 0.3in $\{5/16,9/16\}$
\tskip
$\{1/4,3/4\}$ \hskip 0.55in $R_{3/4}$ \hskip 0.6in $\{3/16,15/16\}$ \hskip 0.3in $\{7/16,11/16\}$
\bigskip
\centerline{
\insertRaster 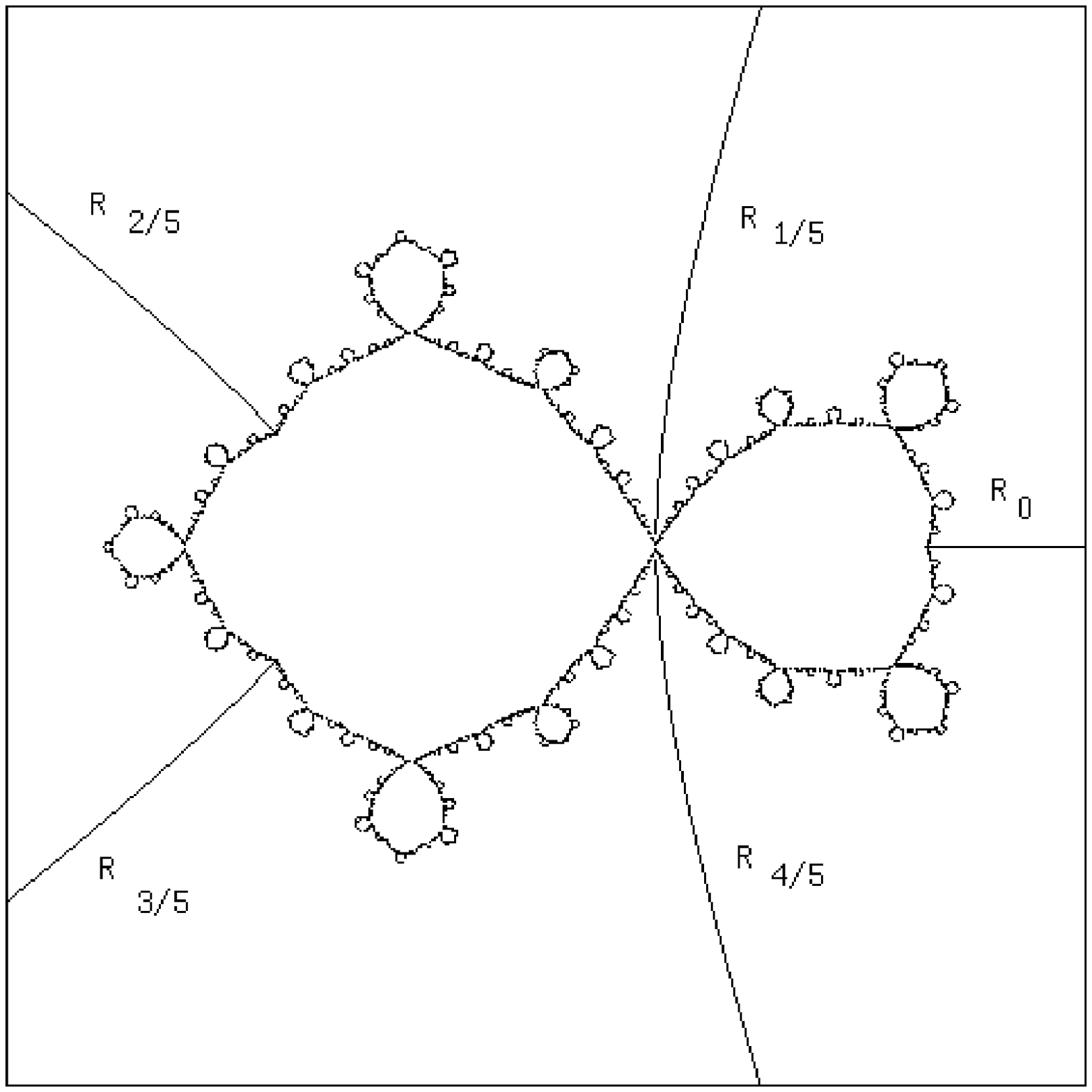 pixels 480 by 480 scaled 337 
\hskip 0.4in 
\insertRaster 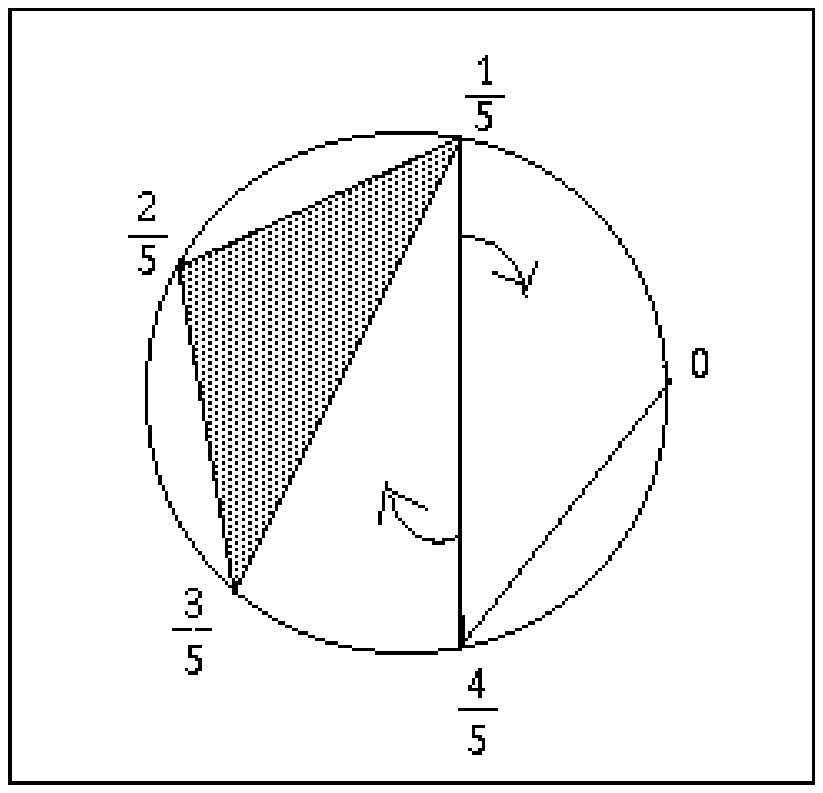 pixels 243 by 234 scaled 600
}
\medskip
\centerline{\it Figure 1.6}
\bigskip 
{\bf 4.7 Badly mixed case.} 
(See Figure 1.6.)
Consider the degree 5 polynomial 
$P(z)=c(z^5+3z^4+3z^3+z^2)$, where $c \approx 4.3582708$. 
It has two Fatou critical components, 
one (on the right) fixed of degree 2, 
and one (on the left) preperiodic of degree 3 (absorbed by the first in one iteration). 
The boundaries of these two Fatou components share a point, which happens to be critical. 
The image of this Julia set critical point is the only fixed point lying in the boundary of the fixed Fatou critical component. 
Only the ray $R_{0}$ lands at this fixed point. 
The rays $R_{1/5},R_{4/5}$ are thus the only rays landing at the Julia set critical point. 
Now, one of these rays ($R_{4/5}$) supports the fixed Fatou component, 
while the other supports the preperiodic one. 
Also $R_{0}$ must have two inverses supporting the fixed Fatou component ($R_{0},R_{4/5}$), 
and three supporting the preperiodic one ($R_{1/5},R_{2/5},R_{3/5}$).  
Thus, the marking is 
${\bf \cal F}=\{\{0,4/5\},\{1/5,2/5,3/5\}\},{\bf \cal J}=\{\{1/5,4/5\}\}$. 
Note that in this case there are arguments that belong to one family and to the other. 
Of course, if this happens, these arguments must be strictly preperiodic. 
\smallskip
By the moment we will take a closer look at condition (c.2) by analyzing this example. 
In this case, conditions (c.1), (c.3)-(c.5) are clearly satisfied. 
To have a degree 5 formal critical portrait, the three sets 
$\{0,4/5\}$, $\{1/5,2/5,3/5\}$, $\{1/5-\epsilon,4/5-\epsilon\}$ 
must be unlinked for $\epsilon>0$ small; which is evidently true. 
\bigskip
{\bf 4.8 Several critical cycles.} 
(See Figure 1.3.) 
Consider the degree 9 polynomial $P \circ P$ where $P$ is as in Example 4.4.
The filled-in Julia set of this polynomial, as well as the external rays remain unchanged (with respect to $P$). 
In this case however, we have two fixed Fatou components each of critical degree 4. 
Each of them absorbs in one iteration another critical component. 
Now each cycle is independent of the other, 
and the choice of markings are independent in the two fixed components. 
Nevertheless, the choice of marking in the fixed Fatou components determines the marking of the critical components they absorb. 
Let us denote by $A,B$ the fixed critical components and by $A',B'$ the critical components they absorb. 
The marking now is ${\bf \cal F}=\{{\cal F}_A,{\cal F}_{A'},{\cal F}_B,{\cal F}_{B'}\},{\bf \cal J}=\emptyset$. 
\bigskip
\hskip 0.05in {\it Component A} \hskip 0.9in ${\cal F}_A$ \hskip 1.6in ${\cal F}_{A'}$ 
\tskip
\hskip 0.3in $R_{3/4}$ \hskip 0.5in $\{3/4,62/72,6/72,14/72\}$
\hskip 0.55in $\{30/72,38/72\}$
\tskip
\hskip 0.3in $R_{7/8}$ \hskip 0.5in $\{7/8,7/72,15/72,55/72\}$
\hskip 0.55in $\{31/72,39/72\}$
\tskip
\hskip 0.3in $R_{1/8}$ \hskip 0.5in $\{1/8,17/72,57/72,65/72\}$
\hskip 0.5in $\{41/72,33/72\}$
\bigskip
\hskip 0.05in {\it Component B} \hskip 0.9in ${\cal F}_B$ \hskip 1.6in ${\cal F}_{B'}$ 
\tskip
\hskip 0.3in $R_{1/4}$ \hskip 0.5in $\{1/4,26/72,42/72,50/72\}$
\hskip 0.5in $\{66/72,2/72\}$
\tskip
\hskip 0.3in $R_{5/8}$ \hskip 0.5in $\{5/8,53/72,21/72,29/72\}$
\hskip 0.5in $\{31/72,39/72\}$
\tskip
\hskip 0.3in $R_{3/8}$ \hskip 0.5in $\{3/8,43/72,51/72,19/72\}$
\hskip 0.5in $\{5/72,69/72\}$
\bigskip
\noindent
This implies that we have 9 possible markings. 
Note that the marking for the components $A,B$ are independent, 
but they uniquely determine the marking for $A',B'$.
\bigskip 
{\bf 4.9} (See Figure 1.7.) 
In our final example we show the importance of working with two separate families ${\cal F},{\cal J}$. 
Consider the sets 
${\cal A}=\{0,{1 \over 3}\}$, ${\cal B}=\{{5 \over 9},{8 \over 9}\}$. 
The polynomial $P(z)=z^3+Az+B$ ($A = 2.25$, $B \approx -0.4330127i$) has marking 
${\cal F}=\{{\cal A},{\cal B}\}$, ${\cal J}=\emptyset$, 
while the polynomial $P(z)=z^3+A'z+B'$ ($A' \approx 2.181104577$, $B' \approx -0.3871686256i$) has marking 
${\cal F}=\{{\cal A}\}$, ${\cal J}=\{{\cal B}\}$.
\bigskip
\centerline{
\insertRaster 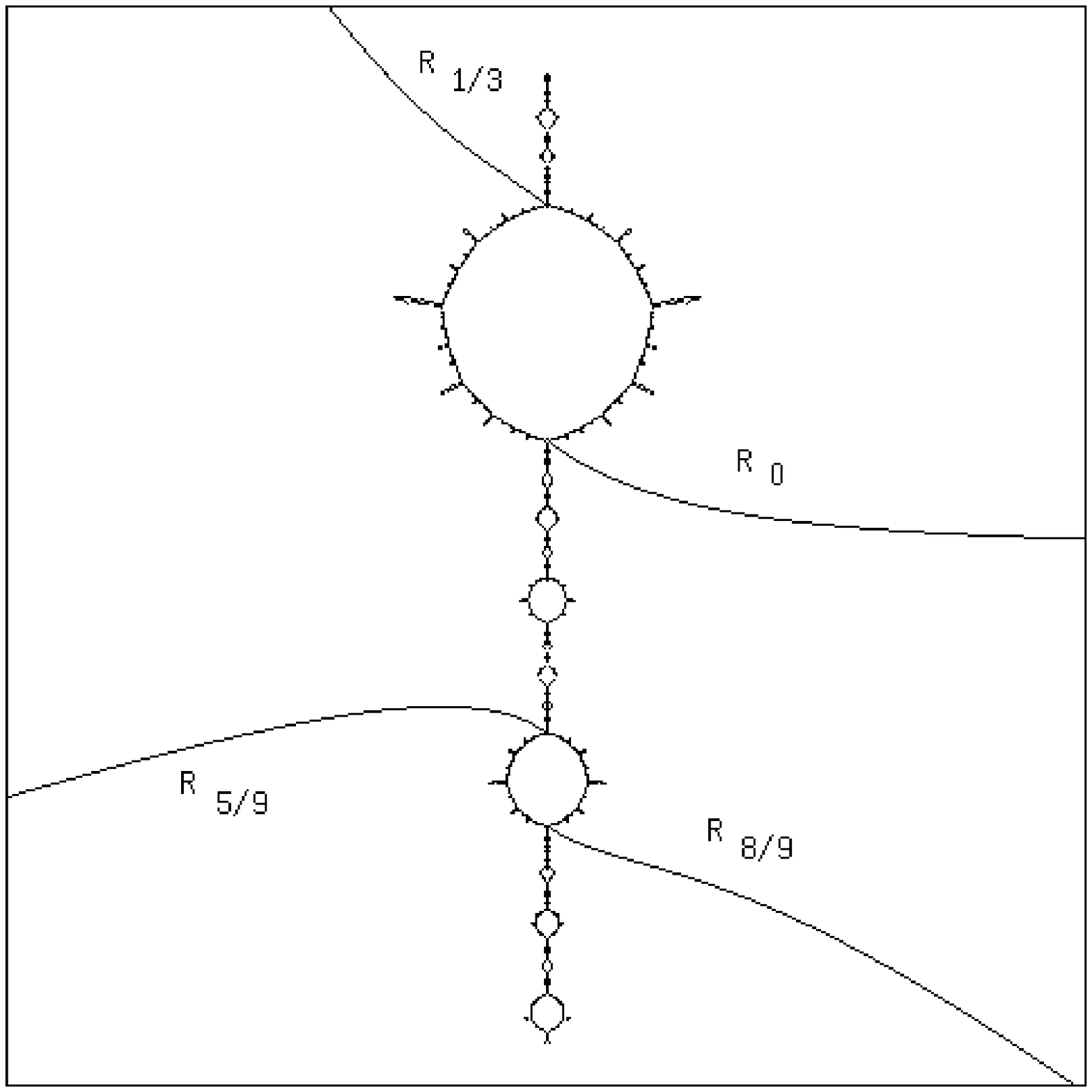 pixels 480 by 480 scaled 337 
\hskip 0.25in 
\insertRaster 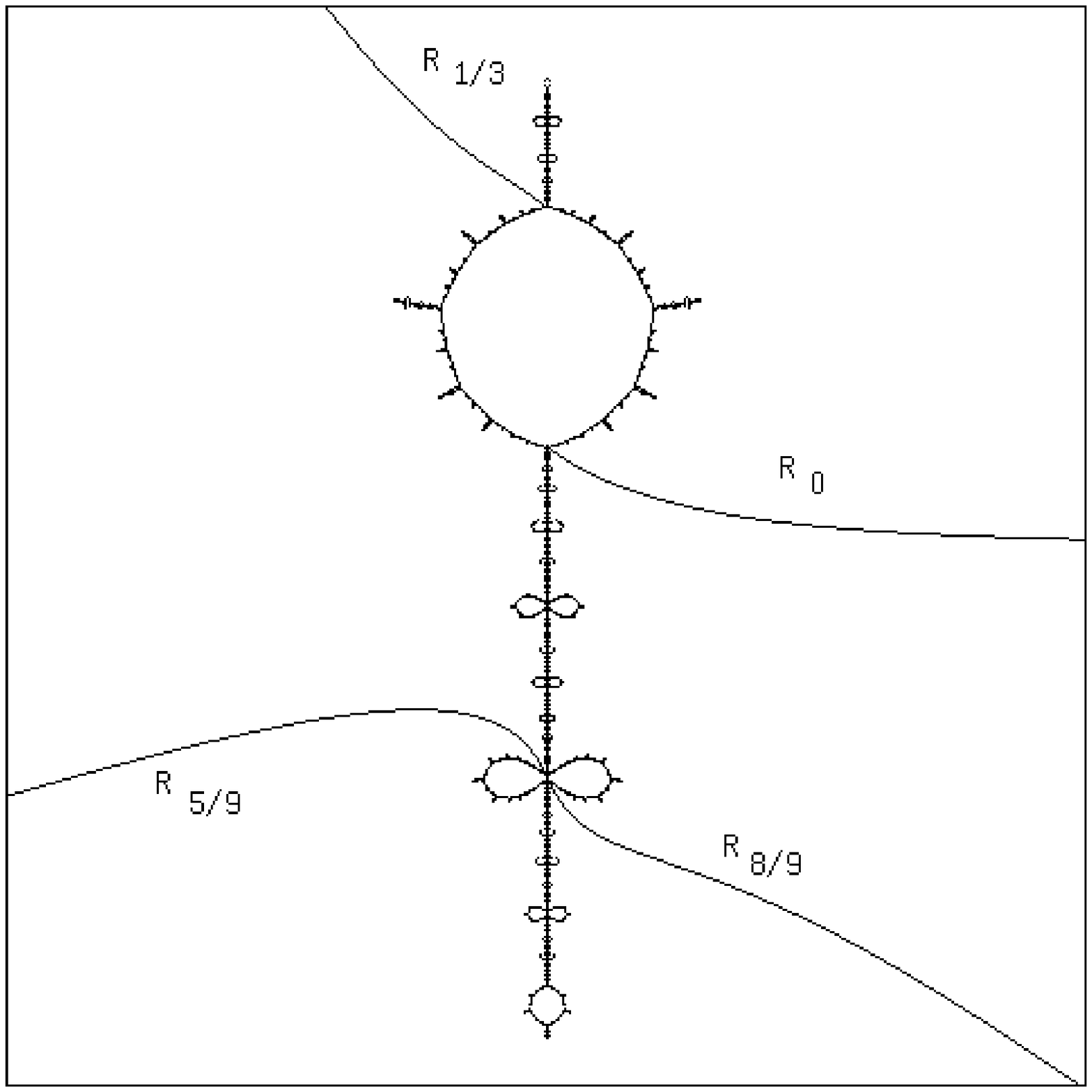 pixels 480 by 480 scaled 337
}
\medskip
\centerline{\it Figure 1.7. Almost the same marking.} 
\vfill\eject

\vfill\eject 
\inchap
\centerline{\bigfont Chapter II }
\centerline{\bigfont Critical Portraits.}
\inchap 
{
In this chapter we isolate the combinatorial properties of a critical portrait $({\cal F},{\cal J})$ as defined in Section I.3, 
and relate them to the dynamics of the respective critically marked polynomial. 
Section 1 deals with the partition in the unit circle determined by this marking. 
We also prove here Lemma I.3.10. 
Section 2 translates to the Julia set the language of Section 1. 
As a consequence we prove that the critical marking defines an admissible critical portrait. 
In Section 3 we prove Proposition I.3.12 
which gives the combinatorial criterion for deciding when two external rays land at the same point. 
Section 4 characterizes the preimages of marked periodic rays landing at that same Fatou component from the combinatorial point of view. 
Almost all the material in this chapter can be found in a weaker formulation in [BFH]. 
The essential novelty here is Section 4, 
which plays a central role in the proof of the realization 
Theorem for Critical Portraits.} 
\insec 
\centerline{\medfont 1. Partitions of the unit circle.}
\insec 
{In this section we fix a formal critical portrait $({\bf {\cal F},{\cal J}})$, 
and study some dynamical properties of the partition determined by these families.} 
\bigskip
Given a formal critical portrait $({\bf {\cal F},{\cal J}})$, 
we defined in Chapter I the partitions ${\cal P}=\{L_1,\dots,L_d\}$ 
and 
${\cal P}^\pm=\{L_1^\pm,\dots,L_d^\pm \}$. 
The first partition omits the arguments in 
${\bf {\cal F}^\cup} \cup {\bf {\cal J}^\cup}$; 
while the other two cover the whole circle ${\bf R/Z}$. 
We also know that each $L_p$ ($L_p^\pm $) is a finite union of open (semiopen) intervals with total length $1/d$ 
(compare Section I.3.6). 
From the dynamical point of view we can say even more. 
\bigskip
{\bf 1.1 Lemma.}
{\it Each $L_p$ is mapped bijectively by $m_d$ onto the complement of a finite set. 
Each $L_p^\pm$ is mapped bijectively by $m_d$ onto the whole unit circle. 
Furthermore these correspondences preserve the circular order}.
\medskip
{\bf Proof.} 
The proof is straightforward and is left to the reader.
\endofproof
\bigskip 
Before the next corollary, 
we recall briefly the standard language for manipulation of symbol sequences.
Let ${\bf S}=(S_0,S_1,\dots)$, where $S_i \in {\cal P}$. 
The {\it shift of ${\bf S}$} is the sequence $\sigma({\bf S})=(S_1,S_2,\dots)$. 
(Formally $\sigma$ is a map from the space of symbol sequences to itself.)
The $i^{th}$ projection $\pi_i$ is the map from symbol sequences to the partition space $\cal P$ defined by $\pi_i({\bf S})=S_i$.
The proof of the following corollary is an easy induction using Lemma 1.1 
and is left to the reader. 
\bigskip
{\bf 1.2 Corollary.} 
{\it Suppose  
$m_d^{\circ n}(\theta)=m_d^{\circ n}(\theta')$ and  $\pi_j(S^+(\theta))=\pi_j(S^+(\theta'))$ for all $j < n$, 
then $\theta=\theta'$. 
(The same is true if we consider left symbol sequences instead.)} 
\endofproof
\bigskip
{\bf Warning.}
Corollary 1.2 is not necessarily true if we compare left with 
right symbol sequences. 
From $S^+(\theta)=S^-(\theta')$ and $m_d(\theta)=m_d(\theta')$, 
we can not infer $\theta=\theta'$. 
For example, in the Ulam-von Neumann map (compare Example I.4.2), 
$S^+(1/4)=S^-(3/4)$, and both arguments become equal after doubling.  
\bigskip
As our partitions are well behaved under iteration, 
it is natural to introduce dynamically defined refinements. 
The fact that these refinements are also unlinked 
allow us derive some basic properties of symbol sequences. 
\bigskip
{\bf 1.3 Definition.}
For $S_0,S_1,\dots \in \cal P$, 
set $U_{S_0,\dots,S_n}=\{\theta \in {\bf R/Z}:\break 
m_d^{\circ i}\theta \in S_i,i=0,...,n\}$. 
The Lebesgue measure of this set is $1/d^{n+1}$ as can be easily verified by induction.
Also set $U_{S_0,S_1,\dots}= \bigcap_{n=0}^{\infty} cl(U_{S_0,\dots,S_n})$. 
This last set being a nested intersection of non empty compact sets, is non empty. 
It is easy to see that if $S(\theta)=(S_0,S_1,S_2,\dots)$ then 
$\theta \in U_{S_0,S_1,S_2,\dots}$. 
It follows that given $S_0,S_1,\dots \in \cal P$, 
there exists an argument which has either left or right symbol sequence $(S_0,S_1,S_2,\dots)$. 
\bigskip
{\bf 1.4 Lemma.} 
{\it For each $n\ge 0$ the family $\{U_{S_0,\dots,S_n}\}$ is unlinked.}
\smallskip
{\bf Proof.} 
This follows by construction and Lemma 1.1.
\endofproof
\bigskip
{\bf 1.5 Lemma.} 
{\it There are only a finite number of arguments which admit a given symbol sequence.}
\smallskip
{\bf Proof.}
Consider the full orbit of both families $\Lambda={\cal O}({\cal F}^\cup) \cup {\cal O}({\cal J}^\cup)$. 
It is enough to prove that 
the number of connected components of $U_{S_0,S_1,\dots,S_n}-\Lambda$ is bounded by a number which depends only on $({\cal F},{\cal J})$. 
We claim that the cardinality $N=\#(\Lambda)$ of $\Lambda$ is the bound we are looking for. 
We prove this by induction. 
For $n=0$ this is clear. 
Now suppose $U_{S_1,S_2,\dots,S_n}-\Lambda=\cup_{i=1}^k I_1$, 
where each $I_\alpha$ is connected and $k \le N$.  
By construction every set $S_0 \cap m_d^{-1}(I_\alpha)$ is completely contained in a component of ${\bf R/Z}-\Lambda$ and therefore is connected. 
The result follows. 
\endofproof 
\bigskip
{\bf 1.6 Lemma.}
{\it Suppose $\theta,\theta'$ have the same periodic left (or right) symbol sequence. 
Then $\theta,\theta'$ are periodic and have the same period.}
\smallskip
{\bf Proof.}
First note that $\theta$ can not be strictly preperiodic. 
For otherwise, eventually it becomes periodic, 
and such periodic argument would have at least two different inverses with the same symbol sequence, in contradiction with Corollary 1.2.  
If $\theta,\theta'$ are periodic of different period, 
we assume without loss of generality that $\theta$ is fixed, but $\theta'$ is not. 
In this case, we have at least three points with the same symbol sequence, 
for which the cyclic order is not preserved under iteration, 
but this is a contradiction to Lemma 1.1. 
Finally, $\theta$ can not be irrational because of Lemma 1.5. 
\endofproof
\bigskip 
{\bf 1.7 Remark.} 
We conclude this section with a trivial remark that will be used later several times. 
If we take $\theta,\theta' \in {\cal J}_k$ 
and $\lambda$ such that $A^-(\lambda)=A^-(\theta)$, 
then by definition $\lambda \in (\theta',\theta]$. 
Analogously,  
if $\theta,\theta' \in {\cal F}_k$ 
and $\lambda$ is such that $A^+(\lambda)=A^+(\theta)$, 
then $\lambda \in [\theta,\theta')$. 
(There is nothing special about ${\cal J}$ or ${\cal F}$ in this formulation; 
but this is the way in which these statements will be used.) 
\eject
\insec 
\centerline{\medfont 2. The induced partitions in the dynamical plane.}
\insec 
{In this section we introduce the induced partition of the Julia set with respect to the given critical marking. 
The main result is that this partition is Markov. 
As a consequence of this, 
we establish that the critical marking of a postcritically finite polynomial 
is in fact an admissible critical portrait, 
establishing in this way Proposition I.3.8.} 
\bigskip
Let $(P,{\bf {\cal F},{\cal J}})$ be critically marked. 
In analogy with the way we constructed a partition ${\cal P}$ of the unit circle 
where only the arguments in ${\cal F}^\cup \cup {\cal J}^\cup$ were omited, 
we will construct a partition of the dynamical plane off the rays with argument in ${\cal J}^\cup$ and extended rays with argument in ${\cal F}^\cup$. 
To simplify this construction we introduce some notation. 
For a set $\Lambda \subset {\bf R/Z}$ we denote by ${\cal R}(\Lambda)$ the set of all external rays with argument in $\Lambda$ and their landing points. 
Also, whenever $\Lambda\subset {\bf R/Z}$ is a set of arguments each of them supporting a Fatou component, 
we denote by ${\cal E}(\Lambda)$ the set of all extended rays with argument in $\Lambda$ and the respective centers of Fatou components. 
\medskip
\centerline{
\insertRaster 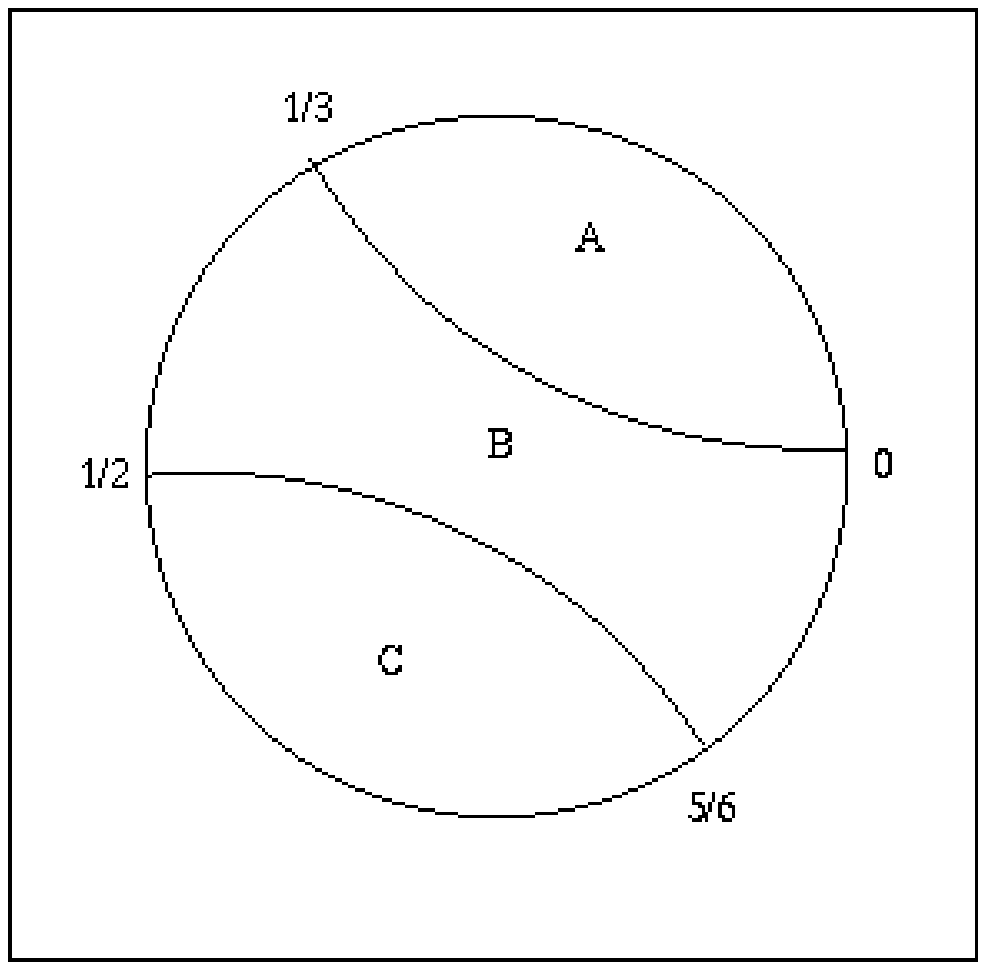 pixels 291 by 286 scaled 500
\hskip 0.4in
\insertRaster 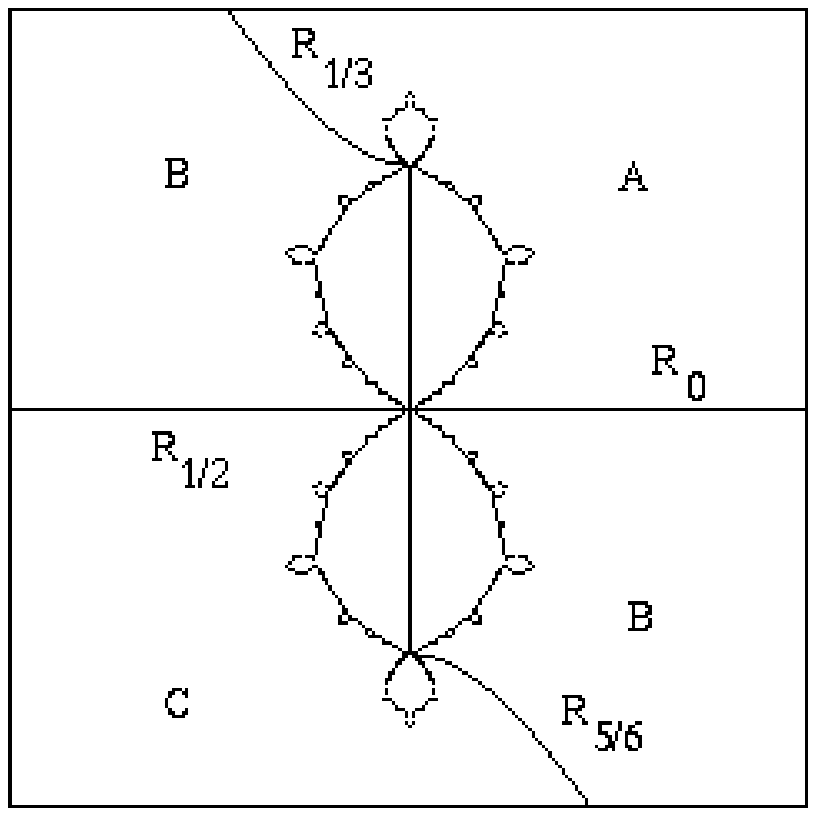 pixels 240 by 240 scaled 650
}
\medskip
\noindent{\it Figure 2.1 The critically marked polynomial $P(z)=z^3+1.5z$ 
with critical portrait $({\cal F}=\{\{0,1/3\},\{1/2,5/6\}\},{\cal J}=\emptyset)$ 
determines a partition of the dynamical plane. 
However the elements of this partition are not necessarily connected open sets. 
Note that 0 and 1/2 share the same left symbol sequence in the circle, 
while the rays $R_0$ and $R_{1/2}$ land at the same point in the dynamical plane.} 
\bigskip
{\bf Definition.} 
We say that two points $z_1,z_2$ in ${\bf C}-{\cal R}({\cal J}^\cup)-
{\cal E}({\cal F}^\cup)$ are {\it ``unlink equivalent"}, 
if they belong to the same connected component of 
${\bf C}-{\cal R}({\cal J}_i)$ and of ${\bf C}-{\cal E}({\cal F}_l)$ for all possible choices of 
${\cal J}_i$ and ${\cal F}_l$ in the marking. 
\bigskip 
Looking at the circle at infinity we immediately derive some properties. 
First, it is easy to see that there are exactly $d$ (=deg $P$) equivalence classes. 
Next, we have that either an external ray is completely contained in an equivalence class, 
or is disjoint from it. 
Furthermore, we have that two rays $R_\theta$ and $R_{\theta'}$ belong to the same equivalence class if and only if 
their arguments $\theta$ and $\theta'$ belong to the same element $S \in {\cal P}$. 
Thus, these equivalence classes are in canonical correspondence with the elements of the partition ${\cal P}$. 
For $S \in {\cal P}$ we denote by ${\cal U}_S$ the corresponding equivalence class in the dynamical plane. 
Each equivalence class is by definition a finite union of unbounded open sets. 
Note that if two arguments belong to the same connected component of some $S \in {\cal P}$, 
then the respective rays will be contained within 
the same connected open region in the dynamical plane. 
\bigskip 
{\bf 2.1 Lemma.} 
{\it 
Each region ${\cal U}_S$ is mapped bijectively by $P$ into the complement of a finite number of rays and extended rays.} 
\endofproof
\bigskip 
{\bf 2.2 Lemma.} 
{\it 
The closure $cl({\cal U}_{S})$ and its restriction to the Julia set $J_{S}=J(P) \cap cl({\cal U}_{S})$ are connected.} 
\endofproof
\medskip 
Both proofs are somehow trivial and are left to the reader 
(compare also the proofs of Lemma 2.3 and Corollary 2.4). 
\bigskip
We can go a step beyond, and take the regions determined by the $n$-fold inverse images of those rays and extended rays. 
Or alternatively we can dynamically define sets ${\cal U}_{S_0,\dots,S_n}$ 
in analogy with $\S$1.3. 
The analogy between this and the definition given in $\S$1.3, is clear: 
by definition, $R_\theta \subset {\cal U}_{S_0,\dots,S_n}$ if and only if $\theta \in U_{S_0,\dots,S_n}$. 
Even if the sets $U_{S_0,\dots,S_n}$ are usually disconnected we have that their closures are not. 
\bigskip
{\bf 2.3 Lemma.} 
{\it 
Let $\gamma:[0,1]\to {\bf C}$ be an arc which crosses neither 
a ray with argument in ${\cal O}(m_d({\cal J}^\cup))$ nor 
an extended ray with argument in ${\cal O}(m_d({\cal F}^\cup))$. 
Suppose further that the image of $\gamma$ is disjoint from the forward orbit of all Fatou critical points. 
If $\gamma$ contains an interior point disjoint from these rays and extended rays, 
then $\gamma$ can be lifted in a unique way within any $cl({\cal U}_S)$, for all $S \in {\cal P}$.} 
\medskip
{\bf Proof.} 
Pick an $S \in {\cal P}$ and 
start the lifting of $\gamma$ at an image point not in the above rays or extended rays. 
Note that the hypothesis guarantees that the lifting can be chosen in such way that it 
never gets into any region ${\cal U}_{S'}$ other than ${\cal U}_S$. 
Uniqueness follows from Lemma 2.1. 
\endofproof
\bigskip
{\bf 2.4 Corollary.} 
{\it 
The closure $cl({\cal U}_{S_0,\dots,S_n})$ and its restriction to the Julia set 
$J_{S_0,\dots,S_n}=J(P) \cap cl({\cal U}_{S_0,\dots,S_n})$ are connected.} 
\medskip
{\bf Proof.} 
Note that if we cut open the plane along all extended rays with argument in ${\cal O}(m_d({\cal F}^\cup))$ and 
remove the forward orbit of all Fatou critical points, 
we are left with a connected set. 
In fact, given a Fatou component $U$, 
there is at most one argument in ${\cal O}(m_d({\cal F}^\cup))$ which supports $U$. 
This follows by construction of critical marking using the hierarchic selection. 
(This is the only place where the hierarchic selection is essentially used in this work!) 
Therefore we can join any two points in the Julia set with a path satisfying the hypothesis of Lemma 2.3. 
The result now follows by induction on $n$. 
\endofproof 
\bigskip 
{\bf Remark.} 
That $J_{S_0,\dots,S_n}$ is connected depends upon the fact that 
the definition of critical marking follows a hierarchic selection. 
Without hierarchic selection for extended supporting rays, 
the statement above is definitely not true. 
\bigskip
At the end, we are mostly interested in the effect of this partition in the Julia set. 
We set $J_{S_0,S_1,\dots}=\bigcap_{n=0}^\infty J_{S_0,\dots,S_n}$ 
Note that because $J(P)$ is locally connected, 
it follows easily that the external ray 
$R_\theta$ lands somewhere in the set 
$J_{S^+(\theta)} \cap J_{S^-(\theta)}$. 
Therefore we should ask if $J_{S(\theta)}$ consists of exactly one point. 
\bigskip 
{\bf 2.5 Lemma.}
{\it For any sequence $(S_0,S_1,\dots)$ the set $J_{S_0,S_1,\dots}$ contains exactly one point.}
\medskip
{\bf Proof} (Compare with [GM, Lemma 4.2]) 
We will make use of the Thurston {\it orbifold metric} associated with $P$. 
Let $M_P$ be the surface with boundary, equal to the disjoint union of all 
$\tilde{\cal U}_S$ defined as $cl({\cal U}_S)$ cut open along all marked rays, extended rays and their forward images, 
and with the orbit of the Fatou critical points removed. 
Define the distance $\rho(z,z')$ between two points of $M_P$ to be the infimum of the lengths 
with respect to the orbifold metric of smooth paths joining $z$ to $z'$ within $M_P$ (or $\infty$ if they belong to different components). 
If $z$ and $z'$ belong to the same subset $J_{S_0,S_1} \subset J(P)$, 
then any path from $P(z)$ to $P(z')$ within $\tilde{\cal U}_{S_1}$
can be lifted back uniquely to a path 
from $z$ to $z'$ within $\tilde{\cal U}_{S_0}$ (compare Lemma 2.3). 
Since the orbifold metric is locally strictly expanding, 
a compactness argument shows that
$$ \rho(P(z),P(z')) \ge c\rho(z,z')$$
for some constant $c>1$, independent of $S_i$ for this $P$. 
Therefore, the inverse map 
$$P^{-1}_{S_0}:J_{S} \mapsto J_{S_0,S}$$
contracts lengths by at least 1/c. 
Hence the iterated inverse images $P^{-1}_{S_0}\circ \dots \circ P^{-1}_{S_n}(J_{S_{n+1}})$ have diameter less than some constant divided by $1/c^n$. 
Taking the limit as $n \to \infty$, we obtain the required unique point.
\endofproof
\bigskip 
{\bf 2.6 Corollary.}
{\it For any sequence $(S_0,S_1,\dots)$ we have 
$P(J_{S_0,S_1,\dots})=J_{S_1,S_2,\dots}$}.
\medskip
{\bf Proof.} 
For some $\theta$, either its left or right symbol sequence $S(\theta)$ equals $(S_0,S_1,\dots)$. 
As the ray $R_\theta$ lands at the unique point contained in $J_{S_0,S_1,\dots}$, the result follows. 
\endofproof
\bigskip
{\bf 2.7 Corollary.} 
{\it If $(S_0,S_1,\dots)$ is a periodic sequence of period m, 
then the unique point in $J_{S_0,S_1,\dots}$ is periodic of period dividing  m.}
\smallskip
{\bf Proof.} 
This follows from Lemma 2.5 and Corollary 2.6. 
In fact, the period is $m$ but this is not a priori obvious, 
this will follow from Proposition 3.6. 
\endofproof
\bigskip
{\bf 2.8 A formal critical portrait not coming from a polynomial.}
Consider the degree 4 formal critical portrait
$$
{\bf \cal J}=\{\{{3 \over 60},{18 \over 60}\},\{{19 \over 60},{34 \over 60}\},\{{1 \over 60},{46 \over 60}\}\},
$$
which does not came from the marking of a polynomial. 
(Compare condition (c.7) in $\S$I.3.7 and Corollary 2.9, here $S^-(19/60)=S^-(46/60)$).
\bigskip
\centerline{
\insertRaster 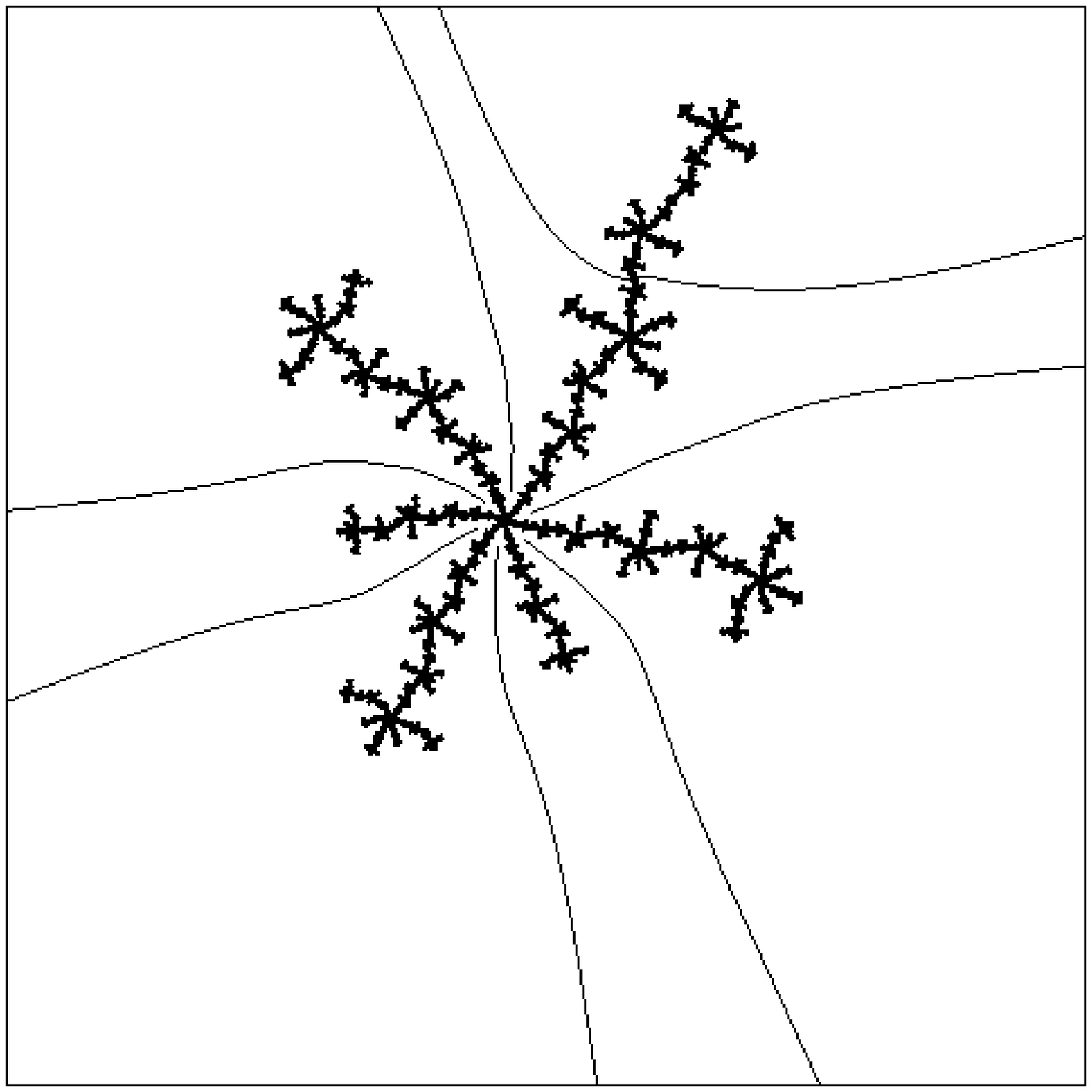 pixels 480 by 480 scaled 360
}

\medskip
\noindent{\it Figure 2.2. 
Julia set of $P(z)=z^4+Az^2+Bz+C$ with the rays 
${1 \over 60}$, ${3 \over 60}$, ${18 \over 60}$, ${19 \over 60}$, ${31 \over 60}$, ${34 \over 60}$, ${46 \over 60}$, ${49 \over 60}$ shown. Here 
\smallskip
\centerline{$A \approx 0.38437710-0.56951210i$}
\centerline{$B \approx 0.30830201+0.03253718i$}
\centerline{$C \approx 0.49119643+0.93292127i$}
}
\medskip
If there is a polynomial $P$ of degree $4$ which realizes this critical portrait, 
there should be critical points $\omega_1 \ne \omega_2$ associated with 
$\{{19 \over 60},{34 \over 60}\}$ and $\{{1 \over 60},{46 \over 60}\}$ respectively. 
But as $S^-(19/60)=S^-(46/60)$, then Lemma 2.5 tell us 
$\omega_1=\omega_2$. 
Thus, the critical points associated with $\{{19 \over 60},{34 \over 60}\},\{{1 \over 60},{46 \over 60}\}$ 
must be actually the same. 
Therefore we do not have three degree 2 critical points, 
but one of degree 3 and the other of degree 2. 
In this case, the rays $R_{4/60}$, and $R_{16/60}$ land at the same fixed point. 
This fixed point has exactly one other preimage, the degree 3 critical point. 
At this critical point the rays $R_{19/60}$, $R_{34/60}$, $R_{49/60}$, $R_{1/60}$, $R_{31/60}$, and $R_{46/60}$ land. 
Therefore, the actual polynomial must have as critical marking either of the following, 
$$
{\cal J}=\{\{{3 \over 60},{18 \over 60}\},\{{19 \over 60},{34 \over 60},{49 \over 60}\}\},
$$
or
$$
{\bf \cal J}=\{\{{3 \over 60},{18 \over 60}\},\{{1 \over 60},{31 \over 60},{46 \over 60}\}\}.
$$
\bigskip
{\bf 2.9 Corollary.}
{\it Let $(P,{\cal F},{\cal J})$ be a critically marked polynomial. 
Suppose $\theta \in {\cal J}_k$ and $\theta' \in {\cal J}_l$. 
If $S^-(m_d^{\circ i}(\theta))= S^-(\theta')$ for some $i \ge 0$, 
then $m_d^{\circ i}(\theta) \in {\cal J}_l$.}
\smallskip
{\bf Proof.} 
It follows from Lemma 2.5 that the rays with argument $m_d^{\circ i}(\theta)$ and $\theta'$ 
land at the same critical point. 
The result then follows from the hierarchic selection of rays. 
(Compare $\S$I.2.) 
\endofproof
\bigskip
{\bf 2.10 Corollary.} 
{\it Let $\gamma \in {\bf {\cal F}^\cup_{per}}$, 
and $\lambda \in {\bf R/Z}$, 
then $\lambda=\gamma$ if and only if\break 
$S^+(\gamma)=S^+(\lambda)$.}
\smallskip
{\bf Proof.} 
Suppose ${\cal F}_k=\{\gamma=\gamma_1,\dots,\gamma_n\}$, 
where the arguments $\gamma_1,\dots,\gamma_n$ are in counterclockwise cyclic order. 
Suppose $\lambda \ne \gamma$ but $S^+(\gamma)=S^+(\lambda)$. 
By Lemma 2.5 the rays $R_\gamma,R_\lambda$ land at the same point. 
As $\lambda$ is periodic by Lemma 1.6, 
it follows that $\lambda \not \in {\cal F}_k$. 
But then, by definition of the right address $A^+(\lambda)$ of $\lambda$, 
it follows that the cyclic order is $\gamma_1,\lambda,\gamma_2,\dots,\gamma_n$ 
(compare Remark 1.7). 
By definition of supporting argument (see $\S$I.1.3), 
the corresponding Fatou component must be in the sector determined by 
$R_{\gamma_1},R_{\lambda}$ (in the counterclockwise sense). 
But this is a contradiction with the fact that $R_{\gamma_2},\dots,R_{\gamma_n}$ also support this component.
\endofproof 
\bigskip
The following now follows from Corollaries 2.9 and 2.10. 
\medskip
{\bf 2.11 Proposition.} 
{\it If $(P,{\bf {\cal F},{\cal J}})$ is a marked polynomial, 
then the pair $({\bf {\cal F},{\cal J}})$ is an admissible critical portrait.}
\endofproof
\vfill\eject
\insec 
\centerline{\medfont 3. Which rays land at the same point?}
\insec 
{We would like to have a combinatorial criterion to decide when two rays land at the same point. 
Two arguments ${\theta,\theta'}$ in the same ${\cal J}_k$ do not have equal (left or right) symbol sequences. 
Nevertheless, the external rays $R_{\theta},R_{\theta'}$ both land at the same critical point. 
In general, all exceptions are a consequence of this fact. 
Furthermore, all the information we need is already contained in left symbol sequences.}
\bigskip
{\bf 3.1 The landing equivalence ($\sim_l$).} 
We recall briefly the definition of the ``landing equivalence" $\sim_l$ between angles, 
introduced in Chapter I 
(compare $\S$I.3.11). 
Let $({\cal F},{\cal J})$ be an admissible critical portrait. 
For $\theta_\alpha,\theta_\beta \in {\cal J}_i \in {\bf \cal J}$ 
we set $S^-(\theta_\alpha) \equiv_i S^-(\theta_\beta)$. 
Then we write $\theta \approx \theta'$ if 
either $S^-(\theta)=S^-(\theta')$ or 
there is an $n \ge 0$ such that $\pi_j(S^-(\theta))=\pi_j(S^-(\theta'))$ for all $j < n$ and 
$\sigma^{n}(S^-(\theta)) \equiv_i \sigma^{n}(S^-(\theta'))$ for some $i$. 
Finally we make this into an equivalence relation by letting 
$\theta \sim_l \theta'$ 
if and only if 
there are arguments $\lambda_0=\theta,\lambda_1,\dots,\lambda_m=\theta'$, 
such that $\lambda_0 \approx \dots \approx \lambda_m$. 
Note that condition $(c.7)$ together with $(c.3)$ guarantee that whenever 
$\theta_i \in {\cal F}_i$ ($i=0,1$); 
then $\theta_0 \sim_l \theta_1$ if and only if ${\cal F}_1={\cal F}_2$. 
\medskip
If the family ${\bf {\cal J}}$ is empty, 
two arguments are equivalent if and only if their left symbol sequences coincide. 
As $S^-(\theta)$ is strictly preperiodic for every argument $\theta$ in the family union ${\bf {\cal J}^\cup}$, 
two periodic or irrational arguments $\theta, \theta'$ are $\sim_l$ equivalent if and only if $S^-(\theta)=S^-(\theta')$. 
Of course, a preperiodic argument would never be 
equivalent to a non preperiodic one. 
\medskip 
By definition, if $\theta \sim_l \theta'$ there is an $m \ge 0$ such that $\sigma^m(S^-(\theta))=\sigma^m(S^-(\theta'))$. 
Also note that whenever $\theta \approx \theta'$ then also $m_d(\theta) \approx m_d(\theta')$.
Therefore the following lemma is trivial. 
\bigskip 
{\bf 3.2 Lemma.} 
{\it If $\theta \sim_l \theta'$ then $m_d(\theta) \sim_l m_d(\theta')$.}
\endofproof
\bigskip 
Now let $(P,{\cal F},{\cal J})$ be a critically marked polynomial. 
We will show now that the $\sim_l$ equivalence classes 
defined from the associated admissible critical portrait effectively characterize the arguments of rays landing at a common point. 
\bigskip
{\bf 3.3 Lemma.}
{\it Suppose $R_\theta,R_{\theta'}$ both land at the same point z. 
If z is non critical then $A^-(\theta)=A^-(\theta')$.} 
\smallskip
{\bf Proof.}
If $z$ is not the landing point of a ray with argument in $\bf {\cal F}^\cup$, 
then it is in the interior of some region ${\cal U}_{S}$. 
Otherwise, 
let $R_{\theta_1},\dots,R_{\theta_k}$ be all rays with argument in ${\cal F}^\cup$ landing at $z$. 
Around $z$ we consider locally segments of these rays 
together with internal rays joining this point $z$ to the center of the $k$ associated Fatou components. 
This configuration divides a neighborhood of $z$ into $2k$ consecutive regions. 
As every other region is contained in ${\cal U}_S$ where 
$S=A^-(\theta_1)=\dots=A^-(\theta_k)$, 
the result follows. 
\endofproof
\bigskip
{\bf 3.4 Lemma.}
{\it Suppose $R_\theta$ lands at a critical point $\omega$, 
then $A^-(\theta)=A^-(\theta')$ for some $\theta' \in {\cal J}_{\omega}$.} 
\smallskip
{\bf Proof.}
The external ray $R_{\theta}$ is contained within some ${\cal U}_{A^-(\theta')}$.
\endofproof
\bigskip
{\bf 3.5 Corollary.}
{\it Suppose $\theta,\theta'$ are such that $A^-(\theta)=A^-(\theta')$. 
Then $R_{\theta},R_{\theta'}$ land at the same point if and only if 
$R_{m_d(\theta)},R_{m_d(\theta')}$ land at the same point.} 
\endofproof
\bigskip
{\bf 3.6 Proposition.} 
{\it Let $(P,{\bf {\cal F},{\cal J}})$ be a marked polynomial. 
Then $R_\theta$ and $R_{\theta'}$ land at the same point if and only if $\theta \sim_l \theta'$.}
\smallskip
{\bf Proof.} 
First suppose that $\theta \sim_l \theta'$. 
If $S^-(\theta) = S^-(\theta')$ then the rays $R_\theta,R_{\theta'}$ land at the same point by Lemma 2.5. 
Otherwise, it is enough to assume $\theta \approx \theta'$. 
In this way, for some $n \ge 0$, $\sigma^n(S^-(\theta)) \equiv_i \sigma^n(S^-(\theta'))$ and $\pi_j(S^-(\theta))=\pi_j(S^-(\theta'))$ for $j < n$. 
By definition there are arguments in ${\cal J}_i$ with symbol sequences $\sigma^n(S^-(\theta))$ and $\sigma^n(S^-(\theta'))$. 
As the rays with these arguments land at the same critical point $\omega_i$, 
the rays $R_{m_d^{\circ n}(\theta)}$ and $R_{m_d^{\circ n}(\theta')}$ also land at $\omega_i$. 
The result follows now from Corollary 3.5. 
\smallskip
Conversely, suppose $R_\theta$ and $R_{\theta'}$ land at the same point $z$. 
There is a minimal $m \ge 0$ such that $P^{\circ m}(z)$ neither is critical nor contains a critical point in its forward orbit. 
We will prove by induction in $m$ that $\theta \sim_l \theta'$. 
Let $P^{\circ n}(z)$ be non critical for all $n \ge 0$ (this is the case $m=0$). 
For all $n \ge 0$, $R_{d^{\circ n}(\theta)},R_{d^{\circ n}(\theta')}$ will be rays landing at the same non critical point. 
In this case the result follows from Lemma 3.3.
Now, let $m_d(\theta) \sim_l m_d(\theta')$ (this is the inductive hypothesis). 
If $z$ is not a critical point we use again Lemma 3.3;  
if $z$ is a critical point we use Lemma 3.4. 
In either case we deduce that $\theta \sim_l \theta'$. 
\endofproof
\bigskip 
{\bf 3.7 Corollary.} 
{\it If $(S_0,S_1,\dots)$ is a periodic sequence of period m, 
then the unique point in $J_{S_0,S_1,\dots}$ has period m.}
\endofproof
\insec 
\centerline{\medfont 4. Which rays support the same Fatou component?}
\insec  
{In general it is impossible to give a combinatorial description of when two arguments support the same Fatou component. 
This because the closure of two Fatou components may share a periodic point 
which is not the landing point of a marked ray. 
In this case, the arguments of all rays landing at such point 
will have the same left and right symbol sequences, 
and thus they are undistinguishable from the combinatoric point of view. 
However, for some cases we will study which rays support some given periodic Fatou component. 
We will only consider rays for which some forward image belongs to the periodic part of the family union $\bf {\cal F}^\cup_{per}$. 
The importance of the combinatorial construction below will become clear in the next chapter. 
In the meanwhile we can tell the reader that in order to apply the theory of ``Levy Cycles" (compare Appendix A), 
we should artificially introduce some preperiodic arguments for every periodic critical point. 
These preperiodic arguments are what we call in this section ``special arguments".}
\bigskip 
As motivation for the combinatorial construction to follow, 
we consider a critically marked polynomial $(P,{\cal F},{\cal J})$. 
Let $\gamma \in {\cal O}({\bf{\cal F}^\cup_{per}})$, be of period $k$. 
Suppose also that $m_d^{\circ nk}(\lambda)=\gamma$. 
\medskip
{\bf 4.1 Lemma.} 
{\it With the above hypothesis, 
$R_\lambda$ supports the same Fatou component as $R_\gamma$ if and only if, 
for each $i \ge 0$ either
\smallskip
i) $\pi_i(S^+\gamma)=\pi_i(S^+\lambda)$, or 
\tskip
ii) $m_d^{\circ i}(\gamma)$ belongs to some ${\cal F}_\alpha$ and $A^+(m_d^{\circ i}(\lambda))=A^+(\gamma')$ for some $\gamma' \in {\cal F}_\alpha$.} 
\smallskip 
{\bf Proof.} 
The proof is straightforward and is left to the reader. 
\endofproof 
\medskip
This motivates the following definition. 
\bigskip
{\bf 4.2 Special arguments.} 
Let $({\bf {\cal F},{\cal J}})$ be an admissible critical portrait. 
To every $\gamma \in {\cal O}({\bf {\cal F}^\cup_{per}})$ we associate a (periodic) sequence of sets ${\cal T}(\gamma,i)$ as follows. 
First we define ${\cal T}(\gamma,0)$: 
$$
{\cal T}(\gamma,0)=\cases{
\{A^+(\gamma'): \gamma' \in {\cal F}_\alpha\} &if $\gamma \in {\cal F}_\alpha$ for some $\alpha$;\cr 
\{A^+(\gamma)\} &otherwise.\cr 
}
$$
In the general case set ${\cal T}(\gamma,j)={\cal T}({m_d^{\circ j}(\gamma),0)}$. 
\bigskip
{\bf Definition.} 
Let $\gamma \in {\cal O}({\bf {\cal F}^\cup_{per}})$ be of period $k=k(\gamma)$. 
We say that $\lambda$ is a {\it special argument} for $\gamma$, 
if there is an $n \ge 0$ such that 
$\pi_i(S^+(\lambda)) \in {\cal T}(\gamma,i)$ for all $i < nk$ and 
$m_d^{\circ nk}(\lambda)=\gamma$. 
In case both $\theta$ and $\theta'$ are special arguments for $\gamma \in {\cal O}({\bf {\cal F}^\cup_{per}})$ 
we write $\theta \sim_\gamma \theta'$. 
\bigskip
The following establishes an equivalence relation between {\it `special arguments'}.
\bigskip
{\bf 4.3 Lemma.} 
{\it If $\lambda$ is a special argument for both $\gamma,\gamma'$, then $\gamma=\gamma'$.}
\smallskip
{\bf Proof.} 
Let $n$ be a multiple of $k(\gamma)k(\gamma')$ big enough, 
then $S^+(\gamma)=\sigma^n(S^+(\lambda))=S^+(\gamma')$ and the result follows from 
condition (c.6) in the definition of admissible critical portrait and Corollary 1.2. 
\endofproof
\bigskip
{\bf 4.4 Remark.} 
If $\theta \sim_\gamma \theta'$ and $S^+(\theta)=S^+(\theta')$, 
it follows from the definition of $\sim_\gamma$, condition (c.6) and Corollary 1.2 that $\theta = \theta'$.
\bigskip
These relations between special arguments are compatible with $m_d$ in the following sense.
\bigskip
{\bf 4.5 Lemma.} 
{\it If $\lambda_1 \sim_\gamma \lambda_2$ then $m_d(\lambda_1) \sim_{m_d(\gamma)} m_d(\lambda_2)$.}
\smallskip
{\bf Proof.} 
For some high iterate $\gamma=m_d^{\circ k}(\lambda_1)=m_d^{\circ k}(\lambda_2)$. 
Thus 
$m_d(\gamma)=m_d^{\circ k}(m_d(\lambda_1))=m_d^{\circ k}(m_d(\lambda_2))$ 
and the result follows from the definition of $\sim_{m_d(\gamma)}$.
\endofproof
\bigskip
The following proposition, is a technical result needed in the proof of 
the main theorem (Theorem I.3.9). 
Its meaning when translated to the context of $PCF$ polynomials, 
is that inverse images of a (marked) periodic ray supporting that same Fatou component, 
can be found very close to the starting periodic ray 
(this is obvious in the context of dynamics, because we are in the subhyperbolic case). 
\bigskip
{\bf 4.6 Proposition.} 
{\it Let $({\bf {\cal F},{\cal J}})$ be an admissible critical portrait.\break 
If $\gamma \in {\bf {\cal F}^\cup_{per}}$ 
then there exist arbitrary small $\epsilon > 0$ such that $\gamma + \epsilon \sim_\gamma \gamma$.}
\medskip
{\bf Proof.} 
Let ${\bf S}_\gamma=(A^+(\gamma),\dots,A^+(m_d^{\circ k-1}(\gamma)))$ and 
take any $W \in {\cal T}^0_\gamma \times \dots \times {\cal T}^{k-1}_\gamma$ 
different from ${\bf S}_\gamma$. 
We form a sequence $\gamma_n \sim_\gamma \gamma$, 
where $S^+(\gamma_n)={\bf S^n_\gamma}W{\bf \bar S_\gamma}$. 
Take a convergent subsequence to $\lambda$. 
As $S^+(\lambda)={\bf \bar S_\gamma}=S^+(\gamma)$ it follows 
by condition (c.6) that $\lambda=\gamma$. 
Now, for $\epsilon > 0$ small enough, 
$\gamma_n$ can not be of the form $\gamma - \epsilon$ by Remark 1.7, 
therefore it must be of the form $\gamma + \epsilon$. 
\endofproof
\bigskip
In the language of special arguments Lemma 4.1 reads. 
\smallskip
{\bf 4.7 Proposition.} 
{\it Let $(P,{\cal F},{\cal J})$ be a marked polynomial. 
If $\theta$ is a\break 
special argument for $\gamma \in {\bf  {\cal F}^\cup_{per}}$ 
then $R_\theta$ and $R_\gamma$ support the same Fatou component.}
\endofproof
\medskip

\vfill\eject  
\inchap
\centerline{\bigfont Chapter III}
\centerline{\bigfont Realizing Critical Portraits}
\inchap 
{In this Chapter we give the proof of the Realization Theorem for Critical Portraits. 
In Section 1 we prove that the combinatorial data is `compatible' in the sense that it allows us 
to construct a Topological Polynomial. 
The actual construction is carried out in Section 2, 
where we also indicate 
(following [BFH]) that it is essentially unique. 
In Section 3 we prove that every admissible critical portrait 
has associated a unique (up to affine conjugation) polynomial which is 
Thurston equivalent to the topological polynomial so far constructed. 
In Section 4 we show that the isotopies between the `actual' and `topological' polynomials 
can be chosen fixed not only relative to certain `marked' points, but also 
relative to the whole boundary when suitably chosen neighborhoods of Fatou points 
are deleted. 
In Section 5 we complete the proof of the Theorem by assigning the expected critical marking to the associated polynomial. }
\insec
\centerline{\medfont 1. Combinatorial Information of Admissible Critical Portraits.}
\insec
{In this Section we analyze the linkage relations that arise when we consider the full orbit of the families and special arguments together. 
The main result is summarized in Proposition 1.2 and is used in Section 2. 
This fact is easy to believe but its proof is extremely technical.} 
\bigskip 
{\bf 1.1} 
Consider an admissible critical portrait $({\bf {\cal F},{\cal J}})$.
The orbit set ${\cal O}({\cal F}^\cup)$ can be partitioned in a natural way as  
${\bf {\cal F}} \cup \{\{\gamma\}:\gamma \in {\cal O}({\cal F}^\cup) - {\cal F}^\cup\}$. 
In the context of dynamics, 
two elements in the orbit ${\cal O}({\cal F}^\cup)$ belong to the same element of this partition if and only if they support the same Fatou component 
(compare Proposition II.4.7). 
If in addition we consider a finite invariant set of special arguments $\Gamma$ 
(i.e, satisfying $m_d(\Gamma) \subset \Gamma \cup {\cal F}^\cup$), 
we can include an element $\lambda \in \Gamma$ in that same class as $\gamma$, 
whenever $\lambda \sim_\gamma \gamma$.
In this way, we construct a family ${\bf {\cal F}^*}=\{{\cal F}^*_1,\dots,{\cal F}^*_n\}$ 
which is a partition of ${\cal O}({\cal F}^\cup) \cup \Gamma$.
\smallskip
Next, we partition the set ${\cal O}({\cal F}^\cup) \cup {\cal O}({\cal J}^\cup) \cup 
\Gamma \cup\{0\}$ 
into $\sim_l$ equivalence classes to form the family ${\bf {\cal J}^*}=\{{\cal J}^*_1,\dots,{\cal J}^*_m\}$. 
In the $PCF$ context we are grouping all those rays we expect to land at the same point 
(compare Proposition II.3.6). 
Here we are adding the argument $\theta=0$ to simplify things later. 
This will reflect the choice of $R_0$ as a preferred fixed `internal' ray in the basin of attraction of $\infty$. 
(Compare Example 3.7.) 
\smallskip
In the way the pair $({{\cal F}^*,{\cal J}^*})$ was constructed, 
it is clear that if we think in terms of external rays, the proposition below must be true. 
\bigskip
{\bf 1.2 Proposition.} 
{\it Let $({\bf {\cal F},{\cal J}})$ be an admissible critical portrait 
and $\Gamma$ a finite invariant set of special arguments. 
With the notation above, 
${\bf {\cal J}^*}$ is weakly unlinked to ${\bf {\cal F}^*}$ in the right.}
\medskip
The reader can skip the rest of this section without any loss of continuity. 
The proof of the proposition follows immediately from Lemmas 1.3-1.9. 
\smallskip 
{\bf 1.3 Lemma.} 
{\it Suppose $\theta_1 \approx \theta_2$, $\psi_1 \approx \psi_2$ but $\theta_1 \not \sim_l \psi_1$. 
Then $\{\theta_1,\theta_2\}$ and $\{\psi_1,\psi_2\}$ are unlinked.}
\smallskip
{\bf Proof.} 
Suppose this is not the case. 
We assume then that $\{\theta_1,\theta_2\}$ and $\{\psi_1,\psi_2\}$ 
are linked because $\theta_2=\psi_2$ implies $\theta_1 \sim_l \psi_1$. 
As a preliminary remark suppose   
$A^-(\theta_1)=A^-(\theta_2)=A^-(\psi_1)=A^-(\psi_2)$; 
then as the cyclic order of these elements is preserved by $m_d$ 
(compare Lemma II.1.1), 
$\{m_d(\theta_1),m_d(\theta_2)\}$ and $\{m_d(\psi_1),m_d(\psi_2)\}$ are still linked. 
For the proof we distinguish several cases.
\medskip
{\it Case 1:} 
$S^-(\theta_1)=S^-(\theta_2)$ and $S^-(\psi_1)=S^-(\psi_2)$. 
This possibility is easily ruled out using Lemma II.1.4. 
We can say even more. 
If $A^-(\theta_1)=A^-(\theta_2)$ and $A^-(\psi_1)=A^-(\psi_2)$
then by that same lemma we have also 
$A^-(\theta_1)=A^-(\psi_1)$. 
Thus, according to our preliminary remark, it is enough to consider the case when 
$A^-(\theta_1) \ne A^-(\theta_2)$.
\medskip
{\it Case 2:} 
$\theta_1,\theta_2 \in {\cal J}_k$. 
As $\psi_1$ and $\psi_2$ belong to different components of 
${\bf R/Z}-{\cal J}_k$, by definition 
$A^-(\psi_1)\ne A^-(\psi_2)$. 
Thus, also by definition $S^-(\psi_1) \equiv_i S^-(\psi_2)$ for some $i$. 
But then, again by definition, there are $\psi_j' \in {\cal J}_i$ ($j=1,2$), 
each in the same connected component of ${\bf R/Z}-\{\theta_1,\theta_2\}$ as $\psi_j$, 
with $S^-(\psi'_j)=S^-(\psi_j)$. 
But this is a contradiction with the fact that ${\cal J}_k,{\cal J}_i$ are unlinked. 
\medskip
{\it Case 3:} 
$S^-(\theta_1) \equiv_k S^-(\theta_2)$ and $S^-(\psi_1)=S^-(\psi_2)$. 
By definition, there is   
 $\theta'_1 \in {\cal J}_k$ such that 
$S^-(\theta'_1)=S^-(\theta_1)$. 
Now, if $\theta_1$ and $\theta_1'$ belong to different components of ${\bf R/Z}-\{\psi_1,\psi_2\}$ then 
$\{\theta_1,\theta'_1\}$, and $\{\psi_1,\psi_2\}$ are linked and we are in case 1. 
Otherwise, we repeat the same reasoning using now $\theta_2$ and we reach either case 1 or case 2. 
\medskip
{\it Case 4:} 
$S^-(\theta_1) \equiv_k S^-(\theta_2)$ and $S^-(\psi_1) \equiv_j S^-(\psi_2)$. 
We proceed as in case 3 and this is reduced to either case 2 or case 3.  
\endofproof
\bigskip
{\bf 1.4 Corollary.} 
{\it The $\sim_l$ equivalence classes are unlinked.}
\endofproof
\bigskip
{\bf 1.5 Lemma.} 
{\it For any ${\cal F}_k \in {\cal F}$ and any $\sim_l$ equivalence class $\Lambda$, 
$\{\Lambda\}$ is weakly unlinked to $\{{\cal F}_k\}$ in the right.}
\smallskip
{\bf Proof.}
Let $\theta_0 \in \Lambda$ and take $\gamma_1,\gamma_2$ consecutive in ${\cal F}_k$ so that $\theta_0 \in (\gamma_1,\gamma_2]$.   
It is enough to prove that if $\theta_0 \approx \theta_1$ then also $\theta_1 \in (\gamma_1,\gamma_2]$. 
If $A^-(\theta_0)=A^-(\theta_1)$, this follows by definition 
($\theta_0$ and $\theta_1$ by definition belong to the same connected component of ${\bf R/Z}-{\cal F}_k$). 
So suppose that $S^-(\theta_0) \equiv_i S^-(\theta_1)$ with $\theta_1 \not\in (\gamma_1,\gamma_2]$. 
In this case there exist ${\cal J}_i \in {\cal J}$ so that $\theta_0' \in {\cal J}_i \cap (\gamma_1,\gamma_2]$ and 
$\theta_1' \in {\cal J}_i \cap (\gamma_2,\gamma_1]$ with 
$S^-(\theta_j) = S^-(\theta_j')$. 
But this is a contradiction with condition (c.2) in the definition of critical portraits 
(${\cal J}_i$ will not be weakly unlinked to ${\cal F}_k$ in the right). 
\endofproof
\bigskip
{\bf 1.6 Lemma.} 
{\it Let $\psi_1 \sim_\gamma \psi_2$ and $\gamma \not \in {\cal F}_k$, 
then $\{\psi_1,\psi_2\}$ and ${\cal F}_k$ are unlinked.}
\smallskip
{\bf Proof.} 
If $A^+(\psi_1)=A^+(\psi_2)$ this follows by definition and Remark II.1.7. 
Otherwise we must have that $\gamma \in {\cal F}_i$ for some $i \ne k$. 
But then a similar argument as that used in Lemma 1.5 shows that 
${\cal F}_i$ and ${\cal F}_k$ are not unlinked. 
\endofproof
\bigskip
{\bf 1.7 Lemma.} 
{\it Let $\theta_i \sim_{\gamma_i} \psi_i$,  $i=1,2$ with $\gamma_1 \ne \gamma_2$. 
Then $\{\theta_1,\psi_1\}$ and $\{\theta_2,\psi_2\}$ are unlinked.}
\smallskip
{\bf Proof.} 
We will consider right symbol sequences $S^+(\theta_j)$ and $S^+(\psi_j)$. 
Suppose is not the case that they are unlinked. 
Then $\{\theta_1,\psi_1\}$ and $\{\theta_2,\psi_2\}$ are linked because $\theta_2=\psi_2$ will imply $\gamma_1 = \gamma_2$ by Lemma II.4.3. 
As preliminary remarks, 
suppose $A^+(\theta_1)=A^+(\psi_1)=A^+(\theta_2)=A^+(\psi_2)$. 
Then as the cyclic order of these elements is preserved by $m_d$ (compare Lemma II.1.1),
$\{m_d(\theta_1),m_d(\theta_2)\}$ and $\{m_d(\psi_1),m_d(\psi_2)\}$ are linked. 
Furthermore, if $A^+(\theta_1)=A^+(\theta_2)$ and $A^+(\psi_1)=A^+(\psi_2)$, 
by Lemma II.1.4 we must have $A^+(\theta_1)=A^+(\psi_1)$. 
\smallskip
Now, suppose $\theta_1$ is in the same connected component of ${\bf R/Z}-\{\theta_2,\psi_2\}$ as $\gamma_1$ (if not $\psi_1$ will be). 
In this case $\{\theta'_1=\gamma_1,\psi_1\}$, and $\{\theta_2,\psi_2\}$ are linked, 
so we assume $\theta_1=\gamma_1$. 
In an analogous way we may suppose that $\theta_2=\gamma_2$. 
Under this assumption we will prove that for all $j \ge 0$, $\{m_d^{\circ j}(\theta_1),m_d^{\circ j}(\psi_1)\}$ and $\{m_d^{\circ j}(\theta_2),m_d^{\circ j}(\psi_2)\}$ should be linked. 
Of course this is absurd because by definition, for $j$ big enough we have $m_d^{\circ j}(\theta_1)=m_d^{\circ j}(\psi_1)=m_d^{\circ j}(\gamma_1)$. 
\medskip
Suppose that $A^+(\theta_1) \ne A^+(\psi_1)$. 
Then by definition $\theta_1 \in {\cal F}_k$ for some $k$. 
Furthermore, there is $\psi'_1 \in {\cal F}_k$ with $A^+(\psi'_1)=A^+(\psi_1)$. 
It follows from Lemma 1.6 that $\theta_1,\psi'_1 \in {\cal F}_k$ are in the same component of ${\bf R/Z}-\{\theta_2,\psi_2\}$. 
Thus, $\{\psi'_1,\psi_1\}$ and $\{\theta_2,\psi_2\}$ are still linked. 
Note that $m_d(\psi')=m_d(\theta_1)$. 
Also by symmetry we may take $A^+(\theta_2)=A^+(\psi_2)$ 
(note that the property $m_d(\theta_2)=m_d(\gamma_2)$ will not be lost). 
But then by the second preliminary remark $A^+(\theta_1)=A^+(\psi_1)=A^+(\theta_2)=A^+(\psi_2)$, 
and so, by the first $\{m_d(\theta_1)=m_d(\gamma_1),m_d(\psi_1)\}$ and $\{m_d(\theta_2)=m_d(\gamma_2),m_d(\psi_2)\}$ are linked. 
This is the desired contradiction. 
\endofproof
\bigskip
{\bf 1.8 Corollary.}
{\it The family $\{\{\theta:\theta \sim_\gamma \gamma\}:\gamma \in {\cal O}({\bf {\cal F}^\cup_{per}})\}$ is\break 
unlinked.}
\endofproof
\bigskip
{\bf 1.9 Lemma.}
{\it Let $\gamma \in {\cal O}({\bf{\cal F}^\cup_{per}})$ and $\Lambda$ an $\sim_l$ equivalence class. 
Then ${\Lambda}$ is weakly unlinked in the right to any finite subset of $\{\theta:\theta \sim_\gamma \gamma \}$.}
\smallskip
{\bf Proof.} 
Take $\gamma \in {\cal F}_{\gamma} \in {\cal F}$. 
We will prove by induction that any $\sim_l$ equivalence class $\Lambda$, 
is weakly unlinked to 
$\Psi_n(\gamma')=\{\theta \sim_{\gamma'} \gamma': m_d^{\circ n}(\theta) \in  {\cal F}_\gamma \}$ 
(here $\gamma'$ belongs to the same cycle as $\gamma$, and $m_d^{\circ n}(\gamma')=\gamma$). 
The result follows easily. 
For $n=0$, this is Lemma 1.5. 
In general take $\theta_1 \approx \theta_2$ and assume that $\{\theta_1,\theta_2\}$ 
is not weakly unlinked in the right to $\{\psi_1,\psi_2\} \subset \Psi_n(\gamma')$. 
\smallskip
{\it Case 1:} $A^+(\psi_1) \ne A^+(\psi_2)$.
Then by definition $\gamma' \in {\cal F}_k\in {\cal F}$ for some $k$. 
Thus, there are $\psi_i' \in {\cal F}_k$ such that $A^+(\psi_i') = A^+(\psi_i)$, 
and because of Lemma 1.5, it is easy to see that 
$\{\theta_1,\theta_2\}$ is not weakly unlinked in the right to either 
$\{\psi_1,\psi_1'\}$ or to $\{\psi_2,\psi_2'\}$ 
(both being subsets of $\Psi_n(\gamma')$). 
Thus it is enough to consider case 2. 
\smallskip
{\it Case 2:} $A^+(\psi_1) = A^+(\psi_2)$.
In this case we can not have simultaneously $\theta_1=\psi_1$ and $\theta_2=\psi_2$. 
In fact, in this case Lemma II.1.1 would imply that  
$\{m_d(\theta_1),m_d(\theta_2)\}$ is not weakly unlinked in the right to $\{m_d(\psi_1),m_d(\psi_2)\}$ in contradiction with the inductive hypothesis. 
Thus we may suppose that $\theta_1 \in (\psi_1,\psi_2)$ (and $\theta_2 \in (\psi_2,\psi_1]$). 
If $A^-(\theta_1)=A^-(\theta_2)$ it follows from Lemma II.1.4 that for $\epsilon >0$ small enough 
$A^+(\theta_1-\epsilon/d)=A^+(\theta_2-\epsilon/d)=A^+(\psi_1)=A^+(\psi_2)$. 
By Lemma II.1.1 we have then that $\{m_d(\theta_1)-\epsilon,m_d(\theta_2)-\epsilon\}$ and $\{m_d(\psi_1),m_d(\psi_2)\}$ are not unlinked, 
in contradiction with the inductive hypothesis. 
Therefore $A^-(\theta_1) \ne A^-(\theta_2)$, and then by definition 
we must have $S^-(\theta_1) \equiv_i S^-(\theta_2)$. 
But then, using the same reasoning as in the previous lemmas, 
we can assume that $\theta_1,\theta_2 \in {\cal J}_i$. 
But if this is the case, 
we get a contradiction 
because it follows by definition and Remark II.1.7 that 
$A^+(\psi_1) \ne A^+(\psi_2)$. 
\endofproof
\bigskip 
Proposition 1.2 follows now easily from the above lemmas. 
\endofproof
\insec 
\centerline{\medfont 2. Abstract and embedded webs.}
\insec 
{In this section we construct from the combinatorial data a topological polynomial of degree $d$. 
We also study some of its basic properties. 
None of the material presented here is essentially new, 
and can be found in a slightly different formulation in [BFH].} 
\bigskip
{\bf 2.1} 
Let $({\bf {\cal F},{\cal J}})$ be an admissible critical portrait. 
For any finite invariant set of special arguments $\Gamma$, 
we consider the pair $({\bf {\cal F}^*,{\cal J}^*})$ as in Section 1. 
With these families, 
we construct first an abstract topological graph $W({\bf {\cal F}^*},{\bf {\cal J}^*})$ as follows. 
We pick a vertex $v=\infty$, and take as many edges ${\cal R}_{\theta}$ 
incident at $\infty$ as elements $\theta \in {\bf{\cal J}^{*\cup}}$. 
Let $v_\theta$ be the other adjacent vertex to ${\cal R}_\theta$. 
We identify the vertices $v_\theta,v_{\theta'}$ if and only if 
$\theta,\theta' \in {\cal J}_k^*$ 
for some $k$; 
that is, 
if and only if $\theta \sim_l \theta'$. 
(This because we are expecting the rays with arguments $\sim_l$ related to land at the same point.) 
We write this vertex as $v({\cal J}_k^*)$. 
As each ${\cal R}_\theta$ is labeled by an argument $\theta$, 
we call it {\it the web ray of argument $\theta$}. 
By abuse of language we will say that $v_\theta$ 
($=v({\cal J}_k^*)$ whenever $\theta \in {\cal J}_k^*$) 
is the {\it landing point of the web ray ${\cal R}_\theta$.} 
\smallskip
Next, for each subset ${\cal F}_k^* \in {\bf {\cal F}^*}$ we consider a new vertex $\omega({\cal F}_k^*)$. 
We join this vertex to the landing points of ${\cal R}_\gamma$ for all $\gamma \in {\cal F}^*_k$. 
(This because, all those rays are supposed to support the same Fatou component; 
compare Proposition II.4.7). 
In this case the {\it extended web ray} ${\cal E}_{\gamma}$ is the set 
formed by the web ray of argument $\gamma$, 
its landing point, and the edge joining this landing point with the vertex 
$\omega({\cal F}_k^*)$. 
In each set ${\cal F}^*_k \in {\cal F}^{*\cup}$ 
there is a preferred argument $\gamma_k$. 
We call the edge $\ell_{{\cal F}_k^*}$ joining 
$\omega({\cal F}_k^*)$ with $v_\gamma$, 
{\it the preferred internal ray associated with the ``Fatou type" point $\omega({\cal F}_k^*)$}. 
\smallskip 
Note that by construction (compare $\S$1.1), 
the argument $0$ is always present in our construction. 
We say that the web ray ${\cal R}_0$ is {\it the preferred internal ray associated with $v=\infty$}. 
The graph $W({\bf {\cal F}^*},{\bf {\cal J}^*})$ constructed in this way, 
is the {\it abstract web} associated with $({\bf \cal F},{\bf \cal J},\Gamma)$. 
We will denote by ${\bf V}$ the set of vertices of this graph. 
\bigskip
{\bf 2.2 Embedded webs.} 
We consider embeddings in the Riemann Sphere ${\bf \hat C}$ 
of this abstract web $W=W({\bf {\cal F}^*,{\cal J}^*})$. 
An embedding such that the cyclic order of the web rays corresponds to 
the cyclic order of the labeling by arguments 
can always be constructed because of Proposition 1.2. 
We can always assume that the respective points at $\infty$ correspond. 
Any such embedding is an {\it embedded web}. 
We still call the image of edges incident at $``\infty"$ {\it web rays}. 
Unless strictly necessary we will not distinguish between an embedding and its image. 
\bigskip
\eject
{\bf 2.3 Web maps.} 
The following two properties follow immediately from the construction of $({\bf {\cal F}^*,{\cal J}^*})$ 
and Lemmas II.3.2 and II.4.5.
\medskip
{\it 
If $\theta,\theta' \in {\cal J}^*_k$, 
there is a unique ${\cal J}^*_{f(k)}$, such that $m_d(\theta),m_d(\theta') \in {\cal J}^*_{f(k)}$.
} 
\tskip
{\it
If $\gamma,\gamma' \in {\cal F}^*_k$, 
there is a unique ${\cal F}^*_{f(k)}$, such that $m_d(\gamma),m_d(\gamma') \in {\cal F}^*_{f(k)}$.
} 
\medskip
These two conditions allow us to define a map $f$ between the set vertices of 
the web $W({\bf {\cal F}^*,{\cal J}^*})$ (also define $f(\infty)=\infty$). 
We can extend this map to a map of the whole graph 
$W({\bf {\cal F}^*,{\cal J}^*})$ as follows. 
For any edge which is a web ray ${\cal R}_\theta$, 
define $f|_{{\cal R}_\theta}$ as an homeomorphism between this edge and the web ray ${\cal R}_{m_d(\theta)}$. 
Otherwise, if 
$\ell$ with adjacent vertices $v_1,v_2$ is not a web ray, 
define $f|_\ell$ as an homeomorphism between this edge and the unique edge with adjacent vertices $f(v_1),f(v_2)$. 
\smallskip 
Note that the above construction determines intrinsically the concept of {\it periodic and preperiodic edges} in the web. 
Also note that preferred internal rays map to preferred internal rays. 
\medskip 
Next, we consider an embedding $\phi:W=W({\bf {\cal F}^*,{\cal J}^*}) \to {\bf \hat C}$.  
Any web map $f$ 
induces a map $\hat f$ of ${\cal W}=\phi(W)$ to itself by the formula 
$$\hat f(z)=\phi(f(\phi^{-1}(z))).$$
By a {\it regular extension of $\hat f$} will be meant 
any extension of $\hat f$ 
which is a degree $d$ orientation preserving branch map of the extended complex plane. 
Keeping track of the embedded vertices this extension is essentially unique. 
\bigskip
{\bf 2.4 Theorem.} 
{\it Let $\phi_1,\phi_2$ be two embeddings of the abstract web  
$W=W({\bf {\cal F}^*,{\cal J}^*})$. 
Let $\hat f_i:{\bf \hat C} \to {\bf \hat C}$ ($i=1,2$) be regular extensions of the web maps. 
Then $(\hat f_1,\phi_1({\bf V}))$ and $(\hat f_2,\phi_2({\bf V}))$ are Thurston equivalent as topological maps (compare Appendix A). 
\smallskip
In fact, 
there are homeomorphisms $\psi_\alpha,\psi_\beta:{\bf \hat C} \to {\bf \hat C}$
isotopic relative to $\phi_1({\bf V})$ so that 
\tskip
i) For every vertex $v \in {\bf V}$, 
$\psi_\alpha(\phi_1(v))=\psi_\beta(\phi_1(v))=\phi_2(v)$. 
\tskip
ii) The diagram 
$$
\matrix{ &\hbox{ }&{\bf \hat C} & \psi_\beta\atop\longrightarrow
& {\bf \hat C} &\hbox{ } \cr
&\hat{f_1} &\Big\downarrow && \Big\downarrow &\hat{f_2} \cr
&\hbox{ }&{\bf \hat C} & \psi_\alpha\atop\longrightarrow 
& {\bf \hat C} &\hbox{ }\cr
}
$$ 
is commutative.}
\smallskip
{\bf Proof.} It is not difficult and can be found in [BFH, Theorem 6.8].
\endofproof
\bigskip
{\bf 2.5 Lifting Webs.}
Suppose ${\cal W}=\phi(W({\bf {\cal J}^*,{\cal F}^*}))$ is an embedded web. 
Given this embedding, 
we fix a regular extension $\hat f:{\bf \hat C} \to {\bf \hat C}$ of the web map. 
If ${\cal W}'$ is another embedded web isotopic to ${\cal W}$ relative to the set $\phi({\bf V})$, 
then $\hat f$ uniquely determines an embedded web ${\cal W}''\subset \hat f^{-1}({\cal W}')$ which is also isotopic to ${\cal W}$ relative to $\phi({\bf V})$, 
as the following construction shows. 
\smallskip 
It is convenient first to define 
``the web ray of argument $0$" in ${\cal W}''$. 
For this we need the following remark. 
\smallskip
{\it 
Let $\theta \ne 0$ belong to ${\cal J}^{*\cup}$. 
If $0 \sim_l \theta$, 
then the web rays ${\cal R}_\theta$ and ${\cal R}_0$ in ${\cal W}$
can not be isotopic relative to the set $\phi({\bf V})$.} 
\smallskip
To see this we note that these web rays determine two sectors. 
By construction each of these two sectors contains all web rays with arguments in $(0,\theta)$ and $(\theta,1)$ respectively. 
Now, by Lemma II.1.6, $\theta$ is of the form $k/(d-1)$, 
so each of the sets 
${\cal J}^{*\cup} \cap (0,\theta)$ and 
${\cal J}^{*\cup} \cap (\theta,1)$ is non empty. 
The result follows easily. 
\smallskip
As a consequence we have that there is a unique edge ${\cal R}_0'$ in ${\cal W}'$ 
which can correspond to ${\cal R}_0$.  
Thus there is a unique `edge' ${\cal R}''_0 \subset {\hat f}^{-1}({\cal R}'_0)$ 
joining $\phi(v_0)$ and $\infty$,  
which is isotopic to ${\cal R}_0$ relative $\phi({\bf V})$. 
This is to be defined as the zero web ray in ${\cal W}''$. 
\smallskip
To construct the web ${\cal W}''$ 
we consider first all edges $\ell \subset {\cal W}$ incident at vertices $v=\phi(v({\cal J}^*_k))$ which are not critical. 
By definition, $\hat f(\ell)$ is also an edge in ${\cal W}$; 
now, there is a unique edge $\ell' \in {\cal W}'$ which is isotopic to 
$\hat f(\ell)$ relative to $\phi({\bf V})$. 
As $\hat f$ is locally one to one near $v$, 
starting at $\hat f(v)$, $\ell'$ can be lifted back in a unique way by $\hat f$ to an arc $\ell''$. 
As $\hat f(\ell)$ and $\ell'$ are in particular isotopic relative to the critical values of $\hat f$, 
it follows that $\ell$ and $\ell''$ are isotopic relative to $\phi({\bf V})$. 
\smallskip 
Finally, we consider all edges $\ell$ incident at critical vertices $v=\phi(v({\cal J}^*_k))$. 
Again we repeat the same procedure but keeping in mind that 
the correct indexing for web rays can be found by its relative position respect to the web ray ${\cal R}_0$. 
The adequate choice of inverses can now be easily determined. 
This finishes the construction of ${\cal W}''$. 
By abuse of notation, 
we denote this embedded web ${\cal W}''$ by ${\hat f}^{-1}({\cal W}')$. 
\medskip
Note that we can apply the same construction to the web ${\cal W}''=\hat f^{-1}({\cal W}')$ and so on; 
in this way we can form a sequence of webs 
$$
{\cal W}',{\hat f}^{-1}({\cal W}'),\dots,{\hat f}^{-n}({\cal W}'),\dots
$$ 
all isotopic relative to $\phi({\bf V})$. 
\insec 
\centerline{\medfont 3. There are no Levy cycles.}
\insec 
{In this Section we will prove that any admissible critical portrait is 
`naturally' associated to a unique polynomial $P$ (see Corollary 3.6). 
The natural way to proceed is to construct from the family $({\bf {\cal F}^*,{\cal J}^*})$ with $\Gamma=\emptyset$ a web map $\hat f$. 
The next step can be (as in [BFH]) to prove that any regular extension has no Thurston's obstruction by proving there are no Levy cycles. 
This fact is by no means obvious. 
In fact, it is easier to prove this fact for maps $\hat f'$ associated to a bigger family $({\bf {\cal F}'^*,{\cal J}'^*})$ 
with $\Gamma$ suitably chosen. 
Now, as a Levy cycle for the map $\hat f$ will determine a Levy cycle for the map $\hat f'$ 
we can conclude that the former map has no Levy cycles.} 
\bigskip 
We start with some notation and another result borrowed from [BFH] Section 7. 
\medskip
{\bf 3.1 Definition.} 
Let $\cal W$ be an embedded web and $\ell \subset {\cal W}$ an edge. 
A Jordan curve ${\cal C}$ disjoint from $\phi({\bf V})$ is said to 
{\it intersect $\ell$ essentially}, 
if for every ${\cal C}'$ homotopic to ${\cal C}$ in ${\bf \hat C}-\phi({\bf V})$, 
we have that $\ell \cap {\cal C}'$ is non empty. 
\medskip
The following is together with Theorem A.5 a technical result needed for the proof of the main theorem. 
\bigskip
{\bf 3.2 Lemma.} 
{\it Suppose $\hat f$ admits a Levy cycle $\Lambda=\{{\cal C}_1,\dots,{\cal C}_k\}$ (see appendix A). 
Then any ${\cal C}_i$ does not intersect a preperiodic edge $\ell$ of the web in an essential way.}
\smallskip
{\bf Proof.} See [BFH] Lemma 7.7. 
\endofproof
\bigskip
{\bf 3.3 Remark.} 
Using Proposition II.4.6 it is easy to construct a finite set of special arguments $\Gamma$ with the following properties.
\smallskip
{\it i) $m_d(\Gamma) \subset \Gamma \cup {\cal O}({\cal F}^\cup)$.}
\tskip
{\it ii) If $\lambda \in {\bf {\cal F}}^\cup_{\bf per}$, 
and $\lambda'$ is the successor (counterclockwise) of $\lambda$ in ${\cal J}^{*\cup}$ 
then $\lambda \sim_{\lambda} \lambda'$.} 
\medskip
In the following lemma we assume that the web and a regular extension 
where constructed with this set of special arguments. 
Here if ${\cal C}$ is a Jordan curve, 
the {\it interior of ${\cal C}$} is defined as the bounded component of ${\bf \hat C}-{\cal C}$.
\bigskip   
{\bf 3.4 Lemma.} 
{\it Let ${\cal C}$ be a Jordan curve disjoint from $\phi({\bf V})$. 
Suppose further that ${\cal C}$ has the following properties, 
\tskip
a) All vertices in $\phi({\bf V})$ which belong to the interior of ${\cal C}$ are periodic and do not belong to a critical cycle. 
\tskip
b) ${\cal C}$ does not intersect essentially any preperiodic edge $\ell$. 
\smallskip
Under theses hypothesis, 
if $v_\theta,v_{\theta'} \in \phi({\bf V})$ 
(corresponding to the landing point of the web rays ${\cal R}_{\theta},{\cal R}_{\theta'}$ respectively) 
belong to the interior of ${\cal C}$, then $A^-(\theta)=A^-(\theta')$.} 
\smallskip
{\bf Proof.} 
Suppose $v_\theta,v_{\theta'}$ are in the interior of ${\cal C}$ 
(and therefore $\theta,\theta'$ are periodic). 
Let $\gamma,\gamma' \in {\cal J}_k$ for some $k$. 
The rays ${\cal R}_{\gamma}$ and ${\cal R}_{\gamma'}$ 
divide the plane in two regions. 
If $v_\theta,v_{\theta'}$ do not belong to the same region, then 
${\cal C}$ will cut either ${\cal R}_{\gamma}$ or ${\cal R}_{\gamma'}$ in an essential way. 
Thus, $\theta,\theta'$ belong to the same connected component of ${\bf R/Z}-{\cal J}_k$. 
Now, let $\gamma,\gamma' \in {\cal F}_k$ for some $k$. 
The extended rays ${\cal E}_{\gamma} $ and ${\cal E}_{\gamma'}$ divide the plane in two regions. 
If both $\gamma,\gamma'$ are preperiodic the same argument as above applies, 
and again $\theta,\theta'$ belong to the same connected component of ${\bf R/Z}-\{\gamma,\gamma'\}$. 
Otherwise, suppose that $\gamma$ is periodic (and thus, $\gamma'$ must be preperiodic). 
By hypothesis there is $\epsilon >0$ such that 
$\gamma+\epsilon$ is a special argument for $\gamma$, 
and 
$(\gamma,\gamma+\epsilon) \cap {\cal J}^{*\cup}= \emptyset$. 
Now we apply the same reasoning 
with the extended rays ${\cal E}_{\gamma+\epsilon} $ and ${\cal E}_{\gamma'}$
and thus 
$\theta,\theta'$ belong to the same connected component of ${\bf R/Z}-\{\gamma+\epsilon,\gamma'\}$ (compare Figure 3.1). 
As $\epsilon$ can be chosen arbitrarily small, 
it follows that for $\epsilon >0 $ small enough, 
$\theta-\epsilon$ and $\theta'-\epsilon$ belong to the same connected component of 
${\bf R/Z}-{\cal F}_k$. 
It follows by definition that $A^-(\theta)=A^-(\theta')$. 
\endofproof
\medskip
{\bf 3.5 Proposition.} 
{\it Let $\hat f:{\bf \hat C} \to {\bf \hat C}$ be a regular extension of the web map over 
$({\bf {\cal F}^*,{\cal J}^*})$ for some $\Gamma$. 
Then $\hat f$ admits no Levy cycles.}
\smallskip
{\bf Proof.}
We are going to add points to $\Gamma$ as needed (see the introduction to this section). 
Suppose by contradiction that $\hat f$ has a Levy cycle 
$\{{\cal C}_1,\dots,{\cal C}_k\}$. 
\smallskip
{\it Step 1.} 
As all ``Fatou points" 
(i.e, vertices of the form $\phi(\omega({\cal F}^*_j))$) 
are preperiodic or belong to a critical cycle, 
no such points are in the interior of an element of a Levy cycle (compare Theorem  A.5). 
\smallskip
{\it Step 2.} 
If $\theta \sim_l \theta'$ but $S^-(\theta) \ne S^-(\theta')$ then $\theta$ is preperiodic, 
and so is $v_\theta$.
Thus, $v_\theta$ is not in the interior of a curve in a Levy cycle.
\smallskip
{\it Step 3.} 
If $v_\theta,v_{\theta'}$ are in the interior of an element of a Levy cycle, 
then by Lemma 3.4 
$A^-(\theta)=A^-(\theta')$. 
\smallskip
{\it Step 4.} 
There are no Levy cycles: 
\smallskip
If $v_\theta,v_{\theta'}$ belong to the interior of an element ${\cal C}_1$ of a Levy cycle, then 
there is another element ${\cal C}$ in this Levy cycle such that 
$v_{m_d(\theta)}$ and $v_{m_d(\theta')}$ belong to the interior of ${\cal C}$. 
This immediately implies $S^-(\theta)=S^-(\theta')$ by step 3 and the definition of Levy cycles. 
In this way $v_\theta=v_{\theta'}$ by construction of the Web. 
But this implies there is a unique point in the interior of an element of a Levy cycle, 
and this is a contradiction with the definition of Levy cycles. 
\endofproof
\medskip
{\bf 3.6 Corollary.} 
{\it Let $({\cal F},{\cal J})$ be an admissible critical portrait. 
There is a unique (up to conjugation) polynomial 
$P({\cal F},{\cal J})$ which is Thurston equivalent to $\hat f$. 
Here $\hat f$ is any regular extension of the web map.} 
\endofproof
\bigskip
{\bf 3.7 Example.} 
We are left with the awkward situation of illustrating a result about the impossibility of Levy cycles. 
In order to do this, 
some hypothesis must be violated. 
We have chosen to violate the condition which avoids the existence of Levy cycles, 
namely that $\sim_l$ equivalence classes determine only one point in the Julia set. 
\smallskip
We consider the admissible critical portrait ${\cal F}=\{\{{1 \over 4},{7 \over 12}\},\{{3 \over 4},{1 \over 12}\}\}$ and ${\cal J}=\emptyset$ 
(compare example I.4.4). 
It is easy to check that $S^-({1 \over 4})=S^-({3 \over 4})$ 
(thus expecting the rays $R_{1 \over 4}$ and $R_{3 \over 4}$ to land at the same point in the Julia set). 
We consider also the set of special arguments 
$\Gamma=\{{13 \over 36},{31 \over 36}\}$ 
which satisfies the hypothesis stated in 3.3 
(here ${13 \over 36} \sim_{1 \over 4} {1 \over 4}$ and ${31 \over 36} \sim_{3 \over 4} {3 \over 4}$). 
Thus we have formed
$$
{\cal F}^*=\{\{{1 \over 4},{13 \over 36},{7 \over 12}\},\{{3 \over 4},{31 \over 36},{1 \over 12}\}\}
$$

$$
{\cal J}^*=\{\{0\},\{{1 \over 12}\},\{{1 \over 4},{3 \over 4}\},\{{13 \over 36}\},\{{7 \over 12}\},\{{31 \over 36}\}\}
$$
(recall the meaning of the elements in each family). 
\smallskip
To illustrate Lemma 3.4 (and Proposition 3.5), 
we construct a web ${\cal W}({\cal F}^*,{\cal J}^*)$ without identifying $v_{1 \over 4}$ and $v_{3 \over 4}$. 
We will show how this leads to a Levy cycle (compare Figure 3.1). 
\smallskip
Lemma 3.4 claims that if there is a Levy cycle, 
then arguments of any two $v_{\theta}$, $v_{\theta'}$ in the interior of a constituent element ${\cal C}$ of this cycle 
should have the same left address. 
In our case this means that any such ${\cal C}$ can not cross any solid segment in Figure 3.1 because of Lemma 3.2. 
Thus, the only possibility of a cycle is as shown in Figure 3.1. 
Of course, with the appropriate identification of $v_{1 \over 4}$ and $v_{3 \over 4}$, 
this is impossible. 
\bigskip
\centerline{
\insertRaster 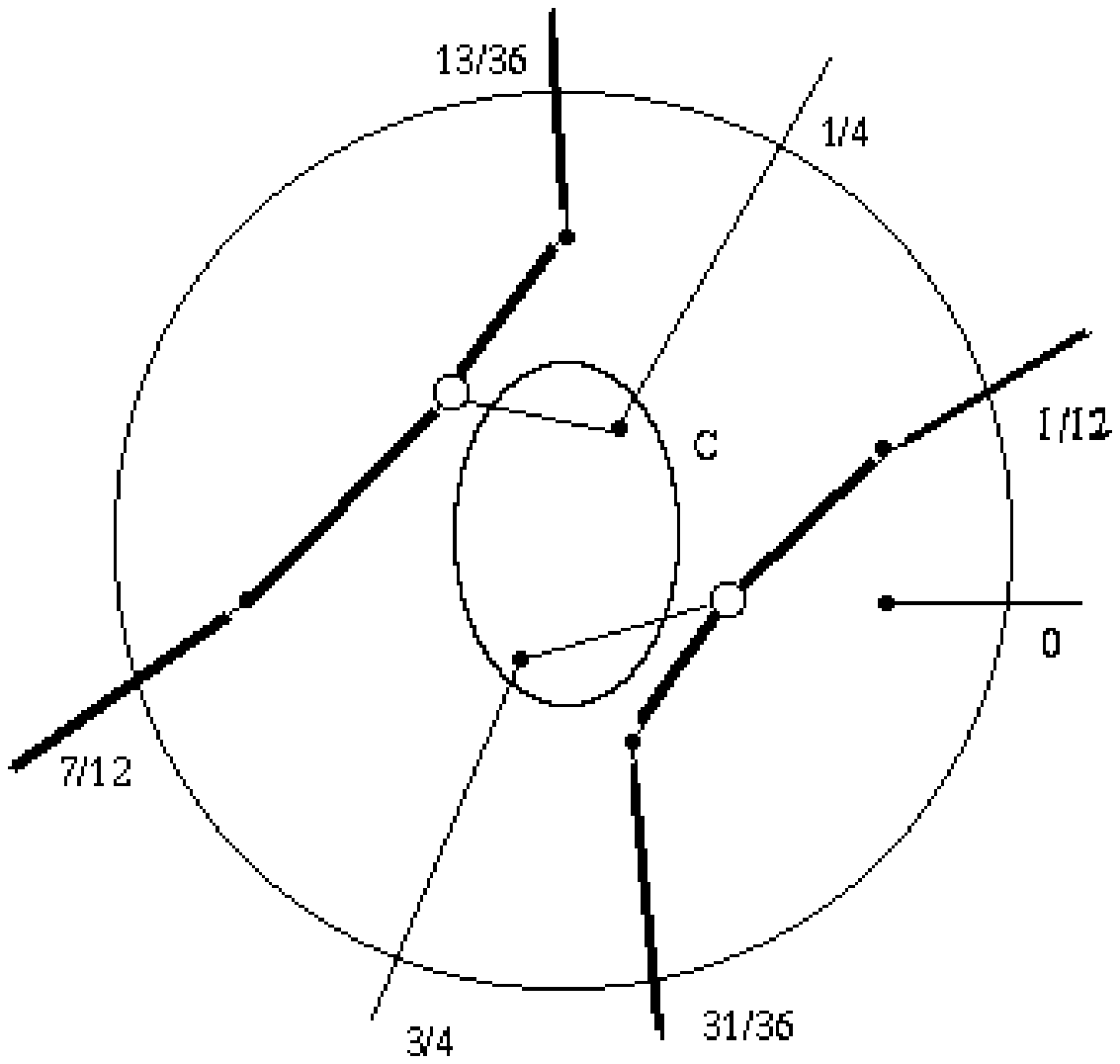 pixels 363 by 325 scaled 360
}
\smallskip
\centerline{\it Figure 3.1}
\insec 
\centerline{\medfont 4. Untwisting the conjugacy.}
\insec 
{Up to this point Corollary 3.6 tells us there is a polynomial (unique up to conjugation) 
associated with the admissible critical portrait $({\cal F},{\cal J})$. 
We must still prove that external and internal rays land at the expected places. 
In other words, 
we have to prove that such post-critically finite polynomial admits the required marking. 
The proof of this fact is not as obvious as it will seem. 
We will consider first a particular example in order to show which difficulties we can still find 
and describe a way to handle them.} 
\bigskip
{\bf 4.1 Example.} 
Consider the admissible critical portrait formed with 
${\cal F}=\{\{0,{1 \over 3},{2 \over 3}\}\}$, ${\cal J}=\emptyset$. 
We first look at the map $f(z)=z^3$ as a {\it `topological polynomial'} 
in the web ${\cal W}({\cal F},{\cal J})$ with {\it vertices} 
${\bf V}=\{0,1,e^{2\pi i \over 3},e^{4\pi i \over 3}\}$ and {\it extended web rays} 
${\cal E}_{k/3}=\{re^{2k\pi i \over 3}:r \in [0,\infty)\}$ for  $k=0,1,2$. 
By Corollary 3.6 this topological polynomial is equivalent to a unique polynomial, 
which will surely be $P(z)=z^3$. 
\smallskip
Consider the homeomorphisms 
$$
\psi_0(r^3e^{2 \pi i\theta})=\cases{
r^3e^{2 \theta \pi i} &if $r \le 3$; \cr 
r^3e^{{2\pi i}[\theta+{3 \over 2}({ln r - ln 3 \over ln 4 - ln 3})]} &if $3 \le r \le 4$; \cr 
r^3e^{{2 \pi i}[\theta+{3 \over 2}]} &if $4 \le r$. \cr
}
$$
\tskip
$$
\psi_1(re^{2 \pi i\theta})=\cases{
re^{2 \theta \pi i} &if $r\le 3$; \cr 
re^{{2 \pi i}[\theta+{1 \over 2}({ln r - ln 3 \over ln 4 - ln 3})]} &if $3 \le r\le 4$; \cr 
re^{{2 \pi i}[\theta+{1 \over 2}]} &if $4 \le r$.\cr
}
$$
Then clearly the following diagram is commutative 
$$
\matrix{ &\hbox{ } &{\bf \hat C} & \psi_1\atop\longrightarrow
& {\bf \hat C} &\hbox{ } \cr
&f &\Big\downarrow && \Big\downarrow  &P. \cr
&\hbox{ } &{\bf \hat C} & \psi_0\atop\longrightarrow
& {\bf \hat C} &\hbox{ } \cr
}\eqno(1)
$$ 
\smallskip
We describe what is happening in the following terms. 
The map $\psi_0$ makes a {\it `Dehn twist'} of 3/2 turns far from $\infty$. 
Thus the `Web'  $\psi_0({\cal W}({\cal F},{\cal J}))$ itself is twisted 3/2 turns. 
By this we mean that when keeping track of the image $\psi_0({\cal R}_0)$ of the web ray ${\cal R}_0$, 
we start as the actual ray $R_0$ for a while, 
then twist in counterclockwise direction until we have completed 3/2 turns, 
and finally continue our way to $\infty$ following the ray $R_{1/2}$! 
Similarly with all other web rays. 
\smallskip
Now, when lifting back the web $\psi_0({\cal W}({\cal F,J}))$ by $P^{-1}$ 
(compare $\S$2.6), 
we see that the resulting embedded web $\psi_1({\cal W}({\cal F,J}))$ has a completely different behavior 
(but they are isotopic). 
The image web ray $\psi_1({\cal R}_0)$ in this case goes for a while in the direction of the actual ray $R_0$, 
then  twists 1/2 turns, 
and finally continues in the direction of the actual ray $R_{1/2}$ to $\infty$. 
\smallskip 
The situation is even worse if we consider successive liftings of the web ray $\psi_0({\cal R}_0)$. 
In these cases, near $\infty$ they will be successively identified with the rays $R_{1 \over 2},R_{1 \over 6},R_{1 \over 18},\dots$. 
Of course, we will prefer to have always near $\infty$ the correct identification. 
In order to describe a possible solution to this dilemma, 
we note that $\psi_0(z)=\psi_1(z)$ for $|z|$ big enough.  
If we remove the set $\{z:|z| > \alpha\}$ for $\alpha$ big enough, 
$\psi_0$ and $\psi_1$ would not be isotopic in this new Riemann surface relative to the boundary 
(they will differ by exactly `one turn' around $\{z:|z|=\alpha\}$. 
This is hardly a surprise because the difference in 1 turn can be easily measured by comparing the embedded web to its lift. 
Now, it is clear that we have not started with the best possible choice of a web.  
Our original web was `twisted' by a given number of turns ($3/2$ in this case); 
when we `lift back' the web, 
this twist will be divided by the degree of the polynomial ($3$ in this case). 
Thus, the `difference in twist' 
(which can always be measured) allows us to state the relation 
$$\hbox{\it twist}-{\hbox{\it twist} \over d}=\hbox{\it difference in twist}.
\eqno(2)$$
Where $d$ is the degree of the polynomial (here $d=3$) and 
{\it difference in twist} is the {\it relative twist of the ray $\psi_1({\cal R}_0)$} in the lifted web respect to the original $\psi_0({\cal R}_0)$. 
In this way, 
equation $(2)$ suggests that any possible odd behavior when lifting webs 
is because of a `Dehn twist' in a neighborhood of Fatou points. 
This is going to be in general the case as we will show below. 
\bigskip 
{\bf 4.2.} 
In the general case, 
we have that starting from the admissible critical portrait $({\cal F},{\cal J})$ 
we can construct a unique up to conjugation polynomial $P$ of degree $d$ 
(which we take here to be monic and centered). 
Also diagram $(1)$ holds. 
Furthermore, 
by replacing $f$ by $\psi_0 \circ f \circ \psi_0^{-1}$ and $\psi_1$ by $\psi_1 \circ \psi_0^{-1}$, 
we may assume without loss of generality that $\psi_0=id$. 
\medskip
For notational convenience we include $\infty$ in the critical set $\Omega(P)$
of the polynomial $P$. 
For each periodic Fatou point $\omega \in \Omega(P)$, 
let $\phi_\omega$ denote a fixed B\"ottcher coordinate associated with $\omega$ ($\infty$ included). 
For $r<1$ define $N_r(w)=\{z \in U(\omega): |\phi_\omega(z)|<r\}$. 
For each strictly preperiodic Fatou point $c \in {\cal O}(\Omega(P))$,  
we inductively define $N_r(c)$ as the connected component of $P^{-1}(N_r(P(c)))$ containing $c$. 
For $X \subset {\cal O}(\Omega(P))$ set  
$N_r(X)=\cup_{c \in X} N_r(c)$. 
\medskip
Now, as there is no topological way to distinguish between the sets 
${\bf \hat C}-{\cal O}(\Omega(P))$ and ${\bf \hat C}-N_r({\cal O}(\Omega(P)))$, 
we can construct an embedded web in ${\bf \hat C}$ and a regular extension $f$ such that the following conditions are satisfied,
\tskip 
i) {\it  $f=P$ in $N_{1/2}({\cal O}(\Omega(P)))$, }
\tskip 
ii) {\it preferred internal web rays are equal to internal preferred rays in $N_{1/2}(\omega)$ if $\omega$ is in a critical cycle,} and  
\tskip
iii) {\it Web edges correspond  to internal rays in $N_{1/2}({\cal O}(\Omega(P)))$.}
\smallskip
Denote by ${\cal W}$ the so constructed web, 
and by ${\bf V}$ be the respective set of vertices 
(there is no further need to write this set as $\phi({\bf V})$). 
Recall we are assuming that $\psi_0$ is the identity in diagram (1). 
Note also that the construction implies that near periodic critical points, 
$\psi_1$ is a rotation in the B\"ottcher coordinate. 
\bigskip 
{\bf 4.3 Untwisting external rays.} 
We consider first what happens near $\infty$ 
(for example, in the set $N_{1/2}(\infty)$). 
As diagram (1) is commutative, we have that for any positive $r \le 1/2$, 
$\psi_0({\cal W}) \cap N_r(\infty)$ is by construction $\infty$ and some segments of actual external rays. 
The portion of the web ray $\psi_1({\cal R}_0) \cap N_r(\infty)$ 
must then be a segment of a ray of the form $R_{j/d}$. 
Furthermore, we can measure the relative twist of $\psi_1({\cal R}_0)$ respect to $\psi_0({\cal R}_0)$ in $\partial N_r(\infty)$ 
(which by construction is a rational number of the form $k/d$).  
Stating this as an equation 
$$\hbox{\it possible twist}-{\hbox{\it possible twist} \over d}=\hbox{\it difference in twist }$$
we have necessarily a rational solution of the form $k/(d-1)$ (same $k$ as above). 
\smallskip
To prove that this `possible twist' is in fact a twist we proceed as follows. 
Take a positive $s<r$ and consider the annulus $N_{r^d}(\infty)-N_{s^d}(\infty)$. 
We modify $\psi_0$ in $N_{r^d}(\infty)$ 
by making a twist of $-{k \over d-1}$ turns inside this annulus. 
This forces us to modify $\psi_1$ in $N_r(\infty)$ 
by a twist of $-{k \over d(d-1)}$ turns 
inside the annulus $N_{r}(\infty)-N_{s}(\infty)$ 
in order to make diagram (1) commutative. 
Clearly there is no problem in doing so because $\psi_0$ is the identity in $N_{r^d}(\infty)$, 
and $\psi_1$ is a rotation in the set $N_r(\infty)$ respect to the B\"ottcher coordinate. 
\smallskip
Formally, we have that in the set ${\bf \hat C}-{\bf V}-N_r(\infty)$, 
$\psi_0$ and $\psi_1$ are not isotopic respect to the boundary because they differ by $k/d$ turns. 
In the annulus $N_r(\infty)-N_{s^d}(\infty)$, the modified $\psi_0,\psi_1$ differ by $-k/d$ turns. 
In this way, 
the modified $\psi_0,\psi_1$ are isotopic relative to the boundary in ${\bf \hat C}-{\bf V}-N_{s^d}(\infty)$. 
Thus, the `difference in twist' between the `new' web rays $\psi_i({\cal R}_0)$ is 0 when measured in $\partial N_{s^d}(\infty)$. 
In particular, if we consider the successive lifting of webs $P^{-n}(\psi_0(W({\bf {\cal F}^*,{\cal J}^*})))$ (compare $\S$2.5), 
all these webs (by construction) will have no difference in twist between the respective lifts of web rays of argument 0. 
We remark that near $\infty$ those web rays are now identified with the ray $R_{-k/d-1}$. 
Also the respective lifting of web rays correspond to bigger and bigger portions of actual rays. 
Of course these rays do not necessarily correspond to the expected ones, 
but they will after conjugation of the polynomial $P$ with $A(z)=e^{-{2k \pi i \over d-1}}z$.
\bigskip
{\bf 4.4 Untwisting periodic preferred internal rays.} 
Our next step is to make the analogue construction in the basin of attraction of finite periodic critical cycles. 
Suppose $\omega_0 \mapsto \omega_1 \mapsto \dots \mapsto \omega_n=\omega_0$ 
is a critical cycle, 
and let $d_i$ be the local degree at $\omega_i$. 
The critical cycle has total degree ${\cal D}=d_0 \times \dots \times d_{n-1}$.
Under the same philosophy as in $\S$4.3 we will like to prove that each coordinate in this cycle was `twisted' by say $x_i$ turns. 
We will denote by $\ell_i$ the preferred internal web edge adjacent to $\omega_i$. 
\smallskip
What we can surely do, 
is to measure the difference in twist when we lift back webs. 
In other words the relative twist of $\psi_1(\ell_i) \subset P^{-1}(\ell_{i+1})$ respect to $\psi_0(\ell_i)$. 
Let this value be $y_i$ 
(which by construction is a rational number with denominator $d_i$). 
If it is true that the coordinates are `twisted', 
then the `possible twist' of $\psi_0(\ell_i)$ is by construction $x_i$; 
while when `lifting back' $\ell_{i+1}$ to get $\psi_1(\ell_i)$, 
its `possible twist' $x_{i+1}$ is divided by $d_i$. 
Thus, if we want to proceed as in $\S$4.3 we must be able to solve the system of equations
$$x_i={x_{i+1} \over d_i} +y_i \hskip 0.6in i=0,\dots,n-1 $$
for $x_i$ rational with denominator ${\cal D}-1$. 
But it is clear that this can be done if we rewrite the system as 
$$
\matrix{
\hbox{$d_0 d_1 \dots d_{n-1}$} &x_0 &= &\hbox{$d_1 \dots d_{n-1}$} &x_1 &+ 
&\hbox{$d_0 d_1 \dots d_{n-1}$} &y_0 \cr
\hbox{$d_1 \dots d_{n-1}$} &x_1 &= &\hbox{$d_2 \dots d_{n-1}$} &x_2 &+ 
&\hbox{$d_1 \dots d_{n-1}$} &y_1 \cr
&\vdots &&& \vdots &&& \vdots \cr
\hbox{$d_{n-2} d_{n-1}$} &x_{n-2} &= &\hbox{$d_{n-1}$} &x_{n-1} &+ 
&\hbox{$d_{n-2} d_{n-1}$} &y_{n-2} \cr
\hbox{$d_{n-1}$} &x_{n-1} &= &&\hbox{$x_0$} &+ 
&\hbox{$d_{n-1}$} &y_{n-1} \cr
       }
$$
\smallskip
With the given solutions $x_0,\dots,x_{n-1}$ 
we proceed to untwist the conjugacy in all neighborhoods of the cycle simultaneously as in $\S$4.3. 
\medskip
{\bf 4.5 Untwisting non periodic Fatou critical components.} 
The last basins that need to be `untwisted' are the ones that correspond to strictly preperiodic Fatou critical points. 
Let $\omega$ be such critical point, and $\omega'=f^{\circ n}(\omega)$ the first critical point in its forward orbit. 
We assume that near $\omega'$ the conjugacy has been already `untwisted'. 
In this case the resulting equation is simply $x_\omega=y_\omega$ so we proceed again as in $\S$4.3. 
\insec 
\centerline{\medfont 5. Proof of Theorem I.3.9.}
\insec 
{\bf 5.1} 
Now we apply successively the construction in $\S$2.5. 
The webs ${\cal W}_n=P^{-n}(\psi_0({\cal W}))$ have edges which coincide with 
the actual internal and external rays in a bigger set after each lifting. 
Given $n$, 
for the web ${\cal W}_n$ 
we consider for each $\bar v$ landing point of ``web rays" and for each edge $\ell$ incident at it, 
the orbifold length of $\ell_n={\bf \hat C} - N_{r^{-d^n}}({\cal O}(\Omega(P)))\cap \ell$. 
For fixed $n$ denote by $\delta_n$ the supremum of such numbers over all possible vertices and edges. 
Note that, as the orbifold metric is strictly expanding for $P$ in ${\bf \hat C} - N_{r}({\cal O}(\Omega(P)))$, 
and each $\ell_n$ is the inverse image of some $\ell'_{n-1}$ we have that $\delta_n \downarrow 0$. 
In this way we have that the respective rays (internal and external) of $P$ 
can be found arbitrarily close to the expected landing points. 
As $J(P)$ is locally connected they actually land there. 
\bigskip
{\bf 5.2} 
To finish the proof of the theorem, 
we only have to prove that the rays $R_\gamma$ associated with a Fatou periodic critical point actually support the respective component. 
But this is trivial if we consider Proposition II.4.6. 
In this case $R_{\gamma},R_{\gamma+\epsilon}$ land in the boundary of the same critical component (compare Proposition II.4.7). 
Thus, in the region determined by the extended rays $\hat R_{\gamma},\hat R_{\gamma+\epsilon}$ 
there is no place for a periodic ray $R_\lambda$ of the same period as $R_\gamma$, if $\epsilon >0$ was chosen small enough. 
This completes the proof of Theorem I.3.9.
\endofproof

\vfill\eject 
\inchap 
\centerline{\bigfont Appendix A} 
\centerline{\bigfont Thurston's Topological Characterization of Rational Maps.}
\inchap
Let $f:S^2 \mapsto S^2$ be an orientation preserving branched 
covering map of the topological sphere. 
The set $\Omega(f)$ of all critical points of $f$ is called the {\it critical set of f}. 
The {\it postcritical set of f} is the set $P({\Omega(f)})=\bigcup^\infty_{n=1} f^{\circ n}\Omega(f)$. 
Whenever the set $P({\Omega(f)})$ is finite, 
we say that $f$ is {\it postcritically finite}.
\medskip
In what follows, we assume always that $f$ is postcritically finite. 
A finite invariant set $M$; 
i.e, $f(M) \subset M$, containing all critical points of $f$ 
is called a {\it marked set.}
In analogy with the previous notation, 
we set $P(M)=\bigcup^\infty_{n=1} f^{\circ n}M$, 
and call it a {\it postmarked set}.
 The elements of $M$ (respectively $P(M)$) are called 
{\it marked points} (respectively {\it postmarked points}). 
We say that $(f,M)$ is a {\it marked branched map}.
\medskip 
Two marked branched maps $(f,M(f))$ and $(g,M(g))$ are 
{\it Thurston\break 
equivalent} if there are homeomorphisms $\phi_1, \phi_2:S^2 \to S^2$, 
isotopic relative to the set $P(M(f))$ 
such that $g \circ \phi_1 = \phi_2 \circ f$, 
and $\phi_1(P(M(f)))=\phi_2(P(M(f)))=P(M(g))$. 
\medskip
We say that a simple closed curve $\gamma \subset S^2-P(M)$ is 
{\it non-peripheral} (for the marked branched map $(f,M)$), 
if each component of $S^2-\gamma$ contains at least two points of $P(M)$. 
A {\it multicurve} $\Gamma=\{\gamma_1,...,\gamma_n\}$ 
is a set of simple, closed, disjoint, 
non-homotopic, non-peripheral curves in $S^2-P(M)$. 
A multicurve $\Gamma$ is {\it stable}, 
if for every $\gamma \in \Gamma$, 
every non-peripheral component of $f^{-1}(\gamma)$ is homotopic (relative to $P(M)$) 
to a curve in $\Gamma$.
\medskip
Let $\gamma_{i,j,\alpha}$ be the components of $f^{-1}(\gamma_j)$ 
homotopic to $\gamma_i$ relative to $P(M)$, 
and $d_{i,j,\alpha}$ be the degree of the 
map $f|_{\gamma_{i,j,\alpha}}:\gamma_{i,j,\alpha} \mapsto \gamma_j$.
We define the $(i,j)$ entry of the {\it Thurston Matrix $f_\Gamma$} as
$$(f_\Gamma)_{i,j}= \sum 1/d_{i,j,\alpha}.$$
Note that by the Perron-Frobenious theorem there is a largest positive 
eigenvalue $\lambda(f_\Gamma)$.
\medskip
There is a smaller function $\nu: P_M \mapsto \{1,2,...,\infty\}$, 
such that for all $x \in f^{-1}(y)$,  
$\nu(y)$ is a multiple of 
$\nu(x)deg_x f$.
We have that the {\it orbifold} $(S^2,P_M,\nu_f)$ is {\it hyperbolic} if its {\it ``Euler characteristic"} satisfies  
$$
2-\sum_{x \in P(M)} (1-1/\nu_f(x)) < 0.
$$ 
Note that we are allowing extensions of the critical and postcritical sets. 
This is because we what to use Thurston's theorem in more 
generality than presented in [DH2] and used in [F] or [BFH]. 
Our marked set is the usual one and maybe a finite number of additional periodic or preperiodic orbits. 
Note that at these additional points, 
the orbifold function has value 1, 
so that the orbifold structure is only determined by the original postcritical set $P(\Omega(f))$. 
\bigskip 
{\bf A.1 Theorem (Thurston's Characterization of Rational Maps).} 
{\it A marked branched map, with hyperbolic orbifold is equivalent to a 
rational function if and only if for any stable multicurve $\Gamma$, 
we have $\lambda(f_\Gamma) < 1$. 
In this case the rational function is unique up to conjugation by a Mobius transformation.}
\medskip
{\bf Proof.} The proof in [DH2] applies without modification.
\endofproof
\bigskip
{\bf A.2 Topological Polynomials.} 
A branched map $f:S^2 \mapsto S^2$ 
is said to be a {\it topological polynomial } if $f^{-1}(\infty)=\infty$.
\medskip
If we are interested only in topological polynomials Thurston's 
theorem is equivalent to the following (see [BFH Theorem 3.2]).
\medskip
{\bf A.3 Theorem.}
{\it A marked topological polynomial (f,M) is equivalent to
 a polynomial if and only if for any stable multicurve $\Gamma$ 
we have $\lambda(f_\Gamma) < 1$. 
In this case, the polynomial is unique up to conjugation by an 
affine transformation.}
\bigskip
{\bf Definition. } A stable multicurve $\Gamma$, with $\lambda(f_\Gamma) 
\ge 1$ is called a {\it Thurston Obstruction (for (f,M)).}
\insec
\centerline{\medfont Levi Cycles}
\insec
Everything here is taken from [BFH] section 4.
\medskip
Let $(f,M)$ be a marked topological polynomial. 
Let $\Gamma$ be a stable multicurve. 
Suppose there exists $\{\gamma_0,...,\gamma_k=\gamma_0\}=
\Lambda \subset \Gamma$ such that for each $i=0,...,k-1$, $\gamma_i$ 
is homotopic relative to $P(M)$ to exactly one component $\gamma'$ of 
$f^{-1}(\gamma_{i+1})$. 
Suppose also that  
$f:\gamma' \mapsto \gamma_{i+1}$ has degree 1. 
Then $\Lambda$ is called a {\it Levy cycle}.
\bigskip
{\bf A.4 Theorem.}
{\it If a marked topological polynomial (f,M) has a 
Thurston obstruction $\Gamma$, then (f,M) has a Levy cycle.} 
\bigskip
{\bf A.5 Theorem.}
{\it The disks of the elements of $\Lambda=
\{\gamma_0,...,\gamma_k=\gamma_0\}$ (i.e, the bounded components of $
S^2-\gamma_i$), contain only cycles of periodic non-critical points of $P(M)$.} 
\bigskip
The last two Theorems together have an interesting interpretation. 
\medskip
{\it For Post-critically finite topological Polynomials, 
only misidentification of periodic points can lead to an obstruction.} 
\vfill\eject
\inchap
\centerline{\bigfont References.}
\inchap
[BFH] B. Bielefeld, Y.Fisher, J. Hubbard, The Classification of Critically Preperiodic Polynomials as Dynamical Systems; 
Journal AMS 5(1992)pp. 721-762.
\smallskip
[DH1] A.Douady and J.Hubbard, \'Etude dynamique des polyn\^omes complexes, part I; 
Publ Math. Orsay 1984-1985. 
\smallskip
[DH2] A.Douady and J.Hubbard, A proof of Thurston's Topological Characterization of Rational Maps; Preprint, Institute Mittag-Leffler 1984. 
\smallskip
[F] Y.Fisher, Thesis; Cornell University, 1989. 
\smallskip
[GM] L.Goldberg and J.Milnor, Fixed Point Portraits; 
Ann. scient. \'Ec. Norm. Sup., $4^e$ s\'erie, t. 26, 1993, pp 51-98 
\smallskip
[M] J.Milnor, Dynamics in one complex variable: Introductory Lectures; 
Preprint \#1990/5 IMS SUNY@StonyBrook.
\smallskip
[P] A.Poirier, On Postcritically Finite Polynomials; 
Thesis,\break 
SUNY@StonyBrook 1993.

\end

\vfill\ejec 
\end